\documentclass[12pt]{article}
\usepackage{curves}
\usepackage{amsmath, amssymb, lamsarrow}
\usepackage{amssymb}

\usepackage{hyperref}
\def\mineappendix{
        \setcounter{section}{1}
        \setcounter{subsection}{0}
        \def\thesection{\Alph{section}}
        \def\sectionap{\@startsection  {section}{1}{\z@}
                        {-3.5ex plus-1ex minus-.2ex} {0ex plus.2ex}
                        {\reset@font\Large\bf  Appendix:  \, }
                        }
        }
\makeatother
\def\Proclaim #1. #2\par{\bigbreak\noindent{\sc#1.\enspace}{\it#2}\par}

\newcommand{\comb}[2]{\left( \begin{array} {c} #1 \\ #2 \end{array} \right)}

\newtheorem{lem}{Lemma}[section]
\newtheorem{cor}[lem]{Corollary}
\newtheorem{thm}[lem]{Theorem}
\newtheorem{pro}[lem]{Proposition}

\newtheorem{defi}[lem]{Definition}
\newtheorem{rem}[lem]{Remark}

\newcommand{\ep}{\epsilon}
\newcommand{\la}{\lambda}

\newenvironment{ack}{\noindent \textbf{Acknowledgments}.}

\title{Schur Q-Polynomials and Kontsevich-Witten Tau Function}

\author{Xiaobo Liu \thanks{Research was partially supported by NSFC grants 11890662 and 11890660.}, Chenglang Yang}

\date{}

\begin{document}

\maketitle

\begin{abstract}
	Using matrix model, Mironov and Morozov recently gave a formula which represents Kontsevich-Witten tau-function as a linear expansion of Schur Q-polynomials.  In this paper, we will show directly that the Q-polynomial expansion in this formula satisfies the Virasoro constraints, and consequently obtain a proof of this formula without using matrix model. We also give a proof for Alexandrov's conjecture that  Kontsevich-Witten tau-function is a hypergeometric tau-function of the BKP hierarchy after re-scaling.
\end{abstract}

\section{Introduction}

Kontsevich-Witten tau-function, denoted by $\tau_{KW}$, is the generating function for intersection numbers of
certain tautological classes on moduli spaces of stable curves (see section \ref{sec:KW} for the precise definition).
It is also the partition function in the two dimensional pure
topological quantum gravity. It was conjectured by Witten and proved by Kontsevich
that this function is a tau function of the KdV hierarchy (cf. \cite{W} and \cite{K}).
This result reveals deep connections between integrable systems
and geometry of moduli spaces of stable curves.
Since Witten's conjecture is equivalent to the Virasoro constraints for $\tau_{KW}$,
this result
 is also the prototype for the Virasoro conjecture for Gromov-Witten invariants
of smooth projective varieties (cf. \cite{EHX} and \cite{CK}). Since KdV hierarchy is a reduction of KP hierarchy, it is
rather natural to consider representing $\tau_{KW}$ in terms of ordinary Schur functions  which provide
polynomial solutions to KP hierarchy (cf. \cite{IZ}).
However it turns out that the coefficients of such representation are very complicated.
In fact, those coefficients are given by determinants of very complicated matrices (cf. \cite{Z} and \cite{BY}).
Using Kontsevich matrix model, Mironov and Morozov recently gave a surprisingly simple and beautiful formula in \cite{MM},
which represents $\tau_{KW}$  as a linear expansion of Schur's Q-polynomials.
In this paper we will give a proof for this formula without using matrix model.

Q-polynomials, also known as Q-functions, were first introduced by Schur in 1911 in the study of projective representations
of symmetric groups. Such functions provide polynomial solutions to BKP hierarchy (cf. \cite{Y} and \cite{KL}).
In some sense, Q-polynomials correspond to square roots of certain ordinary Schur polynomials (see, for example, \cite{O}).
It is very surprising that it turns out to be much simpler to use Q-polynomials rather than ordinary Schur polynomials to expand
Kontsevich-Witten tau-function.

Let $Q_{\la}(\textbf{t})$ be the Q-polynomial associated with a partition
$\la$ where $\textbf{t}=(t_1, t_3, \cdots)$ are independent variables (see section \ref{sec:Q} for the precise definition).
This is a homogeneous polynomial of degree $|\la|$ which only depends on finitely many variables $t_i$ for each fixed partition $\la$.
Define
\begin{equation}  \label{eqn:tauMM}
\tau_{MM}(\textbf{t}):=
\sum_{\lambda\in DP}\left(\frac{\hbar}{16}\right)^{|\lambda|/3} 2^{-l(\la)} \,\, \frac{Q_\lambda(\delta_{k,1})Q_{2\lambda}(\delta_{k,3}/3)}{Q_{2\lambda}(\delta_{k,1})} \,\, Q_\lambda(\textbf{t}),
\end{equation}
where $DP$ denotes the set of strict partitions and $\hbar$ is a formal parameter.
The constant $Q_\lambda(\delta_{k,1})$ is defined to be the value of $Q_\lambda(\textbf{t})$
at the point $t_k=\delta_{k,1}$ for all $k$,
and the constant $Q_{\lambda}(\delta_{k,3}/3)$ is defined similarly.
There is a well known explicit formula for $Q_\lambda(\delta_{k,1})$ which is related
to the hook length formula (see equation \eqref{eqn:hook}). In this paper we will
give an explicit closed formula for $Q_{\lambda}(\delta_{k,3}/3)$ in Theorem \ref{thm:A}.
Hence $\tau_{MM}$ is just a linear combination of
polynomial functions $Q_{\la}(\textbf{t})$. The coefficients of this linear combination are rather simple. They are
explicitly given by Theorem \ref{thm:A} and  equation \eqref{eqn:Q/Q} (see also Remark \ref{rem:A}).

In this paper, we will give a proof for the Virasoro constraints of $\tau_{MM}$ only using properties
of Q-polynomials.
 Define Virasoro operators
\begin{align}
L_m:=\frac{1}{4}\sum_{a+b=2m}\frac{\partial^2}{\partial t_a \partial t_b}
+\frac{1}{2}\sum_{k \geq 1 \atop k \,\, odd}kt_k\frac{\partial}{\partial t_{k+2m}}
-\frac{1}{2\hbar}\frac{\partial}{\partial t_{2m+3}}+\frac{1}{4}t_1t_1\delta_{m,-1}
+\frac{\delta_{m,0}}{16},
      \label{eqn:Lk}
\end{align}
for $m \geq -1$. These operators satisfy the bracket relation
\begin{equation} \label{eqn:VirBr}
 [L_k, L_m] = (k-m) L_{k+m}.
\end{equation}
So they form a half branch of the Virasoro algebra.
\begin{thm}\label{thm:MM=KW}
The function $\tau_{MM}$ satisfies the following Virasoro constraints:
\[ L_m \, \tau_{MM} = 0 \]
for all $m \geq -1$.
\end{thm}
During the proof of this theorem, we also obtain some interesting properties for Q-polynomials. For example,
a multiplication formula for Q-polynomials is given in Lemma \ref{lem:multi pr}.
Formulas for the action of dual Virasoro operators on Q-polynomials are given in Proposition \ref{pro:DL}.
Various non-trivial combinatorial identities for constants $Q_{\lambda}(\delta_{k,3}/3)$ are given in section \ref{sec:Comb}.
Since these constants determine $\tau_{MM}$, these combinatorial identities
 might be useful in the further study of  function $\tau_{MM}$.

Since the Kontsevich-Witten tau function $\tau_{KW}$ also satisfies the same Virasoro constraints which uniquely fix the tau function up to a scalar, an immediate consequence of Theorem \ref{thm:MM=KW} is the following
\begin{cor} \label{cor:KW=MM}
Mironov-Morozov's  formula
\begin{equation} \label{eqn:KW=MM}
\tau_{KW}(\textbf{t})
= \tau_{MM}(\textbf{t})
\end{equation}
 holds.
\end{cor}

Equation \eqref{eqn:KW=MM} was given for Kontsevich matrix model in \cite{MM}. This equation was called the  Mironov-Morozov conjecture in \cite{Alex2020}.  After the first version of this paper has been posted in arXiv:2103.14318, we noticed that a new Remark 3.1 had been added
to the third version of \cite{Alex2020}.
In that remark, it was pointed out (following a communication with Harnad)  that a key step in the derivation of equation \eqref{eqn:KW=MM}, i.e. equation (55) in \cite{MM}, follows from Proposition ($K'$) in \cite{DIZ} and Corolary 3 in \cite{J}.
Actually authors of \cite{DIZ} did not realize that polynomials in their Proposition ($K'$) are Q-polynomials.
This fact was proved in
\cite{J} using a forgotten formula for Hall-Littlewood symmetric functions in \cite{Lit}.

The advantage of our proof of equation \eqref{eqn:KW=MM} using Virasoro constraints is that it does not need matrix model.
Note that Virasoro constraints
exist in a much more general context. For example, the Virasoro conjecture predicts that the generating functions of Gromov-Witten invariants of all smooth projective varieties satisfy Virasoro constraints.  This conjecture has been verified for many special cases
(see, for example, \cite{LT}, \cite{DZ}, \cite{L}, \cite{OP}, \cite{T}). Very few of such geometric models have  matrix model correspondences. Therefore it is of independent interest to have a direct proof for the Virasoro constraints of Q-polynomial series
like $\tau_{MM}$, since this method may be adapted to study similar properties for other geometric models.
For example, in a separate paper \cite{LY} we will  prove
Alexandrov's conjectures posted in \cite{Alex2020} that the Br\'{e}zin-Gross-Witten (abbreviated as BGW) tau function and its generalizations have simple Q-polynomial expansions. It seems that matrix model does not provide much help to understand the Q-polynomial expansion of generalized BGW tau-functions, although these functions were originally defined in matrix model.
The proof in \cite{LY} also used Virasoro constraints, where techniques and some results of this paper
played crucial roles.

The class of hypergeometric tau functions of BKP hierarchy was introduced by Orlov in \cite{O}.
These functions are related to generating functions of spin Hurwitz numbers (cf. \cite{MMN}).
In this paper, we will also give a proof for the following result which was conjectured by Alexandrov
in \cite{Alex2020}:
\begin{thm}\label{cor:KW as BKP}
	The Kontsevich-Witten tau function after a re-scaling, i.e. $\tau_{KW}(\textbf{t}/2)$, is a hypergeometric tau function of the BKP hierarchy.
\end{thm}
This result follows from equation \eqref{eqn:KW=MM} and Corollary \ref{cor:MMhyper}, which states that $\tau_{MM}(\textbf{t}/2)$
is a hypergeometric tau function of BKP hierarchy.
Our proof of Corollary \ref{cor:MMhyper} uses Pfaffian representation of Q-polynomials.
After this work has been completed, we were
informed that Corollary \ref{cor:MMhyper} also follows from the work in \cite{MMNO}
via a completely different approach using fermionic representation of Q-polynomials.

This paper is organized as follows. In section \ref{sec:Pre}, we review definitions and basic properties of Kontsevich-Witten
 tau function and Schur Q-polynomials. In section \ref{sec:dk3}, we compute
 the value of $Q_\la(\delta_{k,3}/3)$ and prove Theorem \ref{cor:KW as BKP}.
 In section \ref{sec:VirMM}, we prove Theorem \ref{thm:MM=KW} and Corollary \ref{cor:KW=MM}. This proof relies on two
 combinatorial identities for the constants $Q_\la(\delta_{k,3}/3)$, which will be proved (together with other identities) in section \ref{sec:Comb}.

\begin{ack}
The authors would like to thank Andrei Mironov and Alexei Morozov for their interests in this work and for bring to our attention several relevant references.
\end{ack}

\section{Preliminaries}
\label{sec:Pre}

\subsection{Kontsevich-Witten tau function}
\label{sec:KW}

Let $\overline{\mathcal{M}}_{g,n}$ be the moduli space of stable genus $g$ curves with $n$ marked points.
Each point in $\overline{\mathcal{M}}_{g,n}$ is represented by a nodal curve $C$ with $n$ distinct smooth
marked points $x_1, \cdots, x_n \in C$, which has finite automorphism group.
For each $i=1, \cdots, n$, let $\mathbb{L}_{i}$ be the line bundle over $\overline{\mathcal{M}}_{g,n}$
whose geometric fiber at a point $(C; x_1, \cdots, x_n) \in \overline{\mathcal{M}}_{g,n}$ is the cotangent
space of $C$ at the marked point $x_i$. Let $\psi_i$ be the first Chern class of $\mathbb{L}_{i}$.
Kontsevich-Witten tau function is the generating function for the intersection numbers
\[ <\tau_{k_1} \tau_{k_2} \cdots \tau_{k_n}>_g \, := \, \int_{\overline{\mathcal{M}}_{g,n}} \psi_1^{k_1} \psi_2^{k_2} \cdots \psi_n^{k_n}, \]
where  $k_1, \cdots, k_n$ are non-negative integers. More precisely, following notations in \cite{Alex2020}, let
\[ \mathcal{F}_{g,n}(\textbf{t}) := \frac{1}{n!} \sum_{k_1, \cdots, k_n \geq 0}  <\tau_{k_1} \tau_{k_2} \cdots \tau_{k_n}>_g
                \prod_{i=1}^n (2k_i+1)!! \, t_{2k_i+1}, \]
where  $\textbf{t}=(t_1, t_3, \cdots)$ are formal variables. Then Kontsevich-Witten tau function
is defined to be
\begin{equation} \label{eqn:KWtau}
\tau_{KW}(\textbf{t}) := \exp \bigg( \sum_{g=0}^{\infty} \sum_{n=0}^{\infty} \hbar^{2g-2+n} \mathcal{F}_{g,n}(\textbf{t})  \bigg),
\end{equation}
where $\hbar$ is a formal parameter.

Note that since the complex dimension of $\overline{\mathcal{M}}_{g,n}$ is $3g-3+n$, $\mathcal{F}_{g,n}(\textbf{t})$
must be a homogeneous polynomial of degree $3(2g-2+n)$ if we assign degree of $t_{k}$ to be $k$ for all $k$.
This implies that in $\log \tau_{KW}(\textbf{t})$, the coefficient of $\hbar^m$ must be a  homogeneous
polynomial of degree $3m$.
By the stability condition, $\overline{\mathcal{M}}_{g,n}$ is non-empty
only if $2g-2+n>0$. Hence only positive integral powers of $\hbar$ can appear in $\log \tau_{KW}(\textbf{t})$.
In particular, at $\textbf{t}=0$, $\tau_{KW}=1$.

It was conjectured by Witten and proved by Kontsevich that $\tau_{KW}$ is a tau function of the KdV hierarchy
(cf. \cite{W} and \cite{K}). Together with the string equation, this property determines the function $\tau_{KW}$.
Let $L_m$ be the Virasoro operators defined by equation \eqref{eqn:Lk}.
Witten's conjecture is equivalent to saying that $\tau_{KW}$ satisfies the following equations:
\begin{equation} \label{eqn:VirKW}
L_m \, \tau_{KW} =0,
\end{equation}
for all $m \geq -1$ (cf. \cite{DVV}, \cite{FKN}, and \cite{KS}).
These equations are known as the {\it Virasoro constraints} for the Kontsevich-Witten tau-function.
In particular, the $L_{-1}$-constraint, i.e. the equation $L_{-1} \, \tau_{KW} =0$, is the string equation.
By the Virasoro bracket relation \eqref{eqn:VirBr}, $L_{-1}$ and $L_2$ generate all operators $L_k$ for
 $k \geq -1$. These two operators are given by
\begin{eqnarray}
	L_{-1} &=& \frac{1}{4}t_1t_1+\frac{1}{2}\sum_{k \geq 1 \atop k \,\, odd} k t_k \frac{\partial}{\partial t_{k-2}}
               -\frac{1}{2\hbar}\frac{\partial}{\partial t_{1}}, \nonumber \\
	L_2 &=& \frac{1}{2}\sum_{k \geq 1 \atop k \,\, odd} k t_k \frac{\partial}{\partial t_{k+4}}
            +\frac{1}{2}\frac{\partial^2}{\partial t_{1}\partial t_{3}}-\frac{1}{2\hbar}\frac{\partial}{\partial t_{7}}.
            \label{eqn:VirOp}
\end{eqnarray}
To prove a function satisfying the Virasoro constraints, we only need to show that it satisfies the
$L_{-1}$ and $L_2$ constraints.

There is also a matrix model, known as Kontsevich model, corresponding to the above geometric description of
Kontsevich-Witten tau-function (c.f. \cite{K}).  The partition function of this matrix model has been shown to
be a tau function of the KdV hierarchy and obeys the Virasoro constraints (c.f. \cite{KMMMZ}, \cite{W91}, and \cite{MMM}).

In \cite{D} and \cite{VV}, it was shown that
recursion relations derived from Virasoro constraints can reduce calculations of all intersection numbers
$<\tau_{k_1} \tau_{k_2} \cdots \tau_{k_n}>_g$ to $<\tau_0^3>_0$ and $<\tau_1>_1$.
In fact, $<\tau_0^3>_0$ and $<\tau_1>_1$ can also be computed from $L_{-1}$ and $L_0$ constraints respectively.
Together with the dimension constraint, this implies that the Virasoro constraints uniquely determine $\tau_{KW}$ up to a constant.
This constant is in fact fixed by the stability condition.
Based on Virasoro constraints,
Alexandrov gave an explicit representation of $\tau_{KW}$ using actions of the cut-and-joint operator, which is a combination of
the slightly modified Virasoro operators with linear coefficients (cf. \cite{Alex11}).
Therefore to prove Corollary \ref{cor:KW=MM}, we only need to show Theorem \ref{thm:MM=KW}, which will be done in section \ref{sec:VirMM}.

\subsection{Schur Q-polynomials}
\label{sec:Q}

Since we have to frequently deal with partitions whose parts might change in the process of calculations,
it will be convenient
not to require partitions to have positive non-increasing parts.
In this paper, a {\it partition} of length $l$ has the form
$\la=(\la_1, \cdots, \la_l)$ with all parts $\la_i$ non-negative integers. We set $l(\la):=l$.
A partition is {\it positive} if all parts are positive.
A partition is {\it weakly positive} if it has at most one part equal to $0$.
A partition $\la$ is {\it strict} if $\la_1 > \la_2 > \cdots > \la_l >0$.
In particular, a strict partition is positive.
The set of all strict partitions is denoted by $DP$.
For convenience, we consider the {\it empty partition} $\emptyset$, i.e. the partition which has no parts,
as a strict partition.

We will follow Macdonald's book \cite{Mac} for the definition of the Schur Q-polynomial $Q_\la$
associated to any partition $\la$.
This is different from the notion used in \cite{MM} and \cite{Alex2020} where the Q-polynomial associated
with $\la$ is equal to $2^{-l(\la)/2} Q_{\la}$.
We consider $Q_{\la}$ as a polynomial of variables $\textbf{t}=(t_1, t_3, t_5, \cdots)$.

More precisely, let $\{q_r(\mathbf{t}) \mid r=0,1,2, \cdots\}$ be the sequence of polynomials whose generating function is given by
\begin{equation} \label{eqn:Ql1}
\sum\limits_{k=0}^\infty q_k(\mathbf{t})z^k=\exp{\bigg(2\sum_{k=0}^{\infty} t_{2k+1} \, z^{2k+1} \bigg)},
\end{equation}
where $z$ is a formal parameter.
If $\lambda=(\lambda_1,...,\lambda_{2m})$ is a weakly positive partition of even length, the associated {\it Schur Q-polynomial}
 is defined by the Pfaffian:
\begin{equation} \label{eqn:Qwp}
Q_\lambda(\textbf{t}):=\text{Pf}(M_{(\lambda_i,\lambda_j)})_{1\leq i,j\leq 2m},
\end{equation}
where
\begin{equation} \label{eqn:Ql2}
M_{(r,s)}:=q_r(\mathbf{t})q_s(\mathbf{t})+2\sum_{i=1}^s(-1)^iq_{r+i}(\mathbf{t})q_{s-i}(\mathbf{t})
\end{equation}
for any pair of non-negative integers $(r,s) \neq (0,0)$ and $M_{(0,0)}:=0$.
It turns out that the matrix $(M_{(\lambda_i,\lambda_j)})_{1\leq i,j\leq 2m}$ is skew symmetric and
the Pfaffian of this matrix is well defined.
For a positive partition $\la$ with odd length,  $Q_{\la}$ is defined to be equal to $Q_{(\la, 0)}$.
Here we have used the notation
$(\lambda, a_1, \cdots, a_k)$ for the partition $(\lambda_1, \cdots, \lambda_l, a_1, \cdots, a_k)$
if $\la=(\la_1, \cdots, \la_l)$ and $a_1, \cdots, a_k$ are non-negative integers.

For convenience, we extend the definition of $Q_{\la}$ to all partitions $\lambda$ by the following rules:
\begin{itemize}
\item If for some $i<l(\la)$, a partition $\tilde{\la}$ is obtained from $\la$ by switching $\la_i$ and $\la_{i+1}$ which are not both equal to $0$, then $Q_{\tilde{\la}} = - Q_{\la}$.
\item $Q_{(\la,0)}= Q_{\la}$ for all $\la$.
\item $Q_{\emptyset}=1$.
\end{itemize}
By properties of Pfaffian, $Q_{\la}$ is skew symmetric with respect to the permutations of $\la$ if $\la$ is weakly positive.
In particular $Q_{\la}=0$ if $\la$ has two equal positive parts.
However $Q_{\la}$ does not change sign when switching two parts of $\la$ which are both $0$.

It is well-known that $Q_\lambda(\textbf{t})$ is a homogeneous polynomial of degree $|\lambda|:=\sum_{i=1}^{l(\lambda)}\lambda_i$,
 where the degree of $t_k$ is assigned to be $k$. Moreover $\{Q_\lambda(\textbf{t}) \mid \lambda\in DP\}$ is a basis of $\mathbb{Q}[t_1,t_3,...]$.
There is a standard inner product $\langle,\rangle$ on the space $\mathbb{Q}[t_1,t_3,...]$ such that
\begin{align}\label{eqn:inner prod}
\langle Q_\lambda,Q_\mu\rangle=2^{l(\lambda)}\delta_{\lambda,\mu}, \text{\ if\ } \lambda, \mu\in DP.
\end{align}
For any operator $f$ on $\mathbb{Q}[t_1,t_3,...]$,  the adjoint operator $f^{\perp}$ is defined by
 \[ \langle f^{\perp} p_1(\textbf{t}), p_2(\textbf{t}) \rangle = \langle p_1(\textbf{t}), f p_2(\textbf{t}) \rangle \]
 for all $p_1, p_2 \in \mathbb{Q}[t_1,t_3,...]$.
In particular, for any positive odd integer $r$, the adjoint of the operator defined by multiplication by $t_r$
is given by
\begin{align}\label{eqn:ad-relation}
t_r^\perp=\frac{1}{2r}\frac{\partial}{\partial t_r}.
\end{align}
 The action of this operator on Q-polynomials can be computed in the following way (cf. \cite{Mac}, p266):
For any strict partition $\lambda$ and positive odd integer $r$, we have
\begin{equation}\label{eqn:dQ}
	\frac{1}{2}\frac{\partial}{\partial t_r} Q_\lambda=\sum_{i=1}^{l(\lambda)}Q_{\lambda-r\epsilon_i},
\end{equation}
where ${\lambda-r\epsilon_i}:=(\lambda_1,...,\lambda_i-r,...,\lambda_l)$. In case $\la_i<r$, $Q_{\lambda-r\epsilon_i}$ should be
understood in the following way:
\begin{defi} \label{def:Q-}
Assume $\la=(\la_1, \cdots, \la_l) \in \mathbb{Z}^l$ with exactly one $\la_i<0$ and all $\la_j \geq 0$ for $j \neq i$.
If there exists $j>i$ such that $\la_j = -\la_{i}$ and $\la_k \neq - \la_i$ for all $k>i$ and $k \neq j$, then define
\[ Q_{\lambda}
	:= (-1)^{j-i-1+\lambda_j} \, 2 \,Q_{\la^{\{i,j\}}}. 	\]
Otherwise we define $Q_{\lambda}:=0$.
\end{defi}
In the above formula, we have used the following notation
\begin{equation}
\lambda^{\{i_1, \cdots, i_n\}}:=(\lambda_1, \cdots, \widehat{\lambda}_{i_1}, \cdots, \widehat{\lambda}_{i_n}, \cdots, \lambda_l)
\end{equation}
for any integers $1 \leq i_1 < \cdots < i_n \leq l$.

\begin{lem}\label{lem:multi pr}
For any strict partition $\lambda$ and positive odd integer $r$, we have
\begin{align}\label{eqn:multi pr}
rt_rQ_\lambda=\sum_{i=1}^{l(\lambda)}Q_{\lambda+r\epsilon_i} +\frac{1}{2}Q_{(\lambda,r)}
+\sum_{k=1}^{r-1}\frac{(-1)^{r-k}}{4}Q_{(\lambda,k,r-k)}.
\end{align}
\end{lem}
{\bf Proof}:
Since $\{Q_\lambda(\textbf{t}) \mid \lambda\in DP \}$ is an orthogonal basis of $\mathbb{Q}[t_1,t_3,...]$,
\begin{equation} \label{eqn:tr}
 rt_rQ_\lambda = \sum_{\mu \in DP} 2^{-l(\mu)} \, \langle rt_rQ_\lambda, \, Q_{\mu} \rangle \, Q_{\mu}.
\end{equation}
By skew symmetry of Q-polynomial, we can replace $\mu$ by any permutation of $\mu$ in the above formula.

By equations \eqref{eqn:ad-relation} and \eqref{eqn:dQ}, we have
\begin{equation} \label{eqn:tradj}
 \langle rt_rQ_\lambda, \, Q_{\mu} \rangle
    =\langle Q_\lambda,\frac{1}{2}\frac{\partial}{\partial t_r}  Q_{\mu} \rangle
    =  \sum_{j=1}^{l(\lambda) } \langle Q_\lambda,  Q_{\mu-r\epsilon_j} \rangle.
\end{equation}
By orthogonality of Q-polynomials, $\langle Q_\lambda,  Q_{\mu-r\epsilon_j} \rangle$ is non-zero only
in the following three cases:

Case (1): $\mu$ is a permutation of $\la+r\epsilon_i$ for some $1 \leq i \leq l(\la)$. In this case, we can replace
$\mu$ by $\la+r\epsilon_i$ in equation \eqref{eqn:tr} and use equation \eqref{eqn:tradj} to compute
$\langle rt_rQ_\lambda,Q_{\lambda+r\epsilon_i} \rangle$.

Note that $\langle Q_{\la}, Q_{\lambda+r\epsilon_i-r\epsilon_j} \rangle \neq 0$ only if
$\lambda+r\epsilon_i-r\epsilon_j$ is a permutation of $\la$. This implies $i=j$ since otherwise $\la_i+r = \la_k$ for
some $k \neq i$,  which is not possible since $\la+r\epsilon_i$ is a permutation of a strict partition.
Hence by equation \eqref{eqn:tradj},
\begin{equation} \label{eqn:trc1}
\langle rt_rQ_\lambda,Q_{\lambda+r\epsilon_i} \rangle
=\langle Q_\lambda,Q_\lambda \rangle
 =2^{l(\lambda)}.
\end{equation}

Case (2): $\mu$ is a permutation of $(\la, k, r-k)$ for some integer $k$ with $0<k<r/2$.
In this case, we can replace $\mu$ by $(\la, k, r-k)$.
Since $(\la, k, r-k)$ is a permutation of a strict partition,
both $k$ and $r-k$ can not be parts of $\la$. This implies
$\langle Q_\lambda, Q_{(\lambda-r\epsilon_j,k,r-k)} \rangle =0$ for $1 \leq j \leq l(\la)$.
Since $Q_{(\la, k, -k)}=0$ by definition, by equation \eqref{eqn:tradj},
\begin{equation} \label{eqn:trc2}
\langle rt_rQ_\lambda,Q_{(\lambda,k,r-k)} \rangle
=\langle Q_\lambda, Q_{(\lambda,k-r,r-k)} \rangle
 =(-1)^{r-k} 2^{l(\lambda)+1}.
\end{equation}

Case (3): $\mu$ is a permutation of $(\la, r)$.
In this case, we can replace $\mu$ by  $(\la, r)$.
Since $(\la, r)$ is a permutation of a strict partition, $r$ can not be a part of $\la$.
Since $\la$ is a strict partition, $0$ is not a part of $\la$.
Hence
$\langle Q_\lambda, Q_{(\lambda-r\epsilon_j,r)}  \rangle = 0$ for all $j \leq l(\la)$.
By equation~\eqref{eqn:tradj},
\begin{equation} \label{eqn:trc3}
\langle rt_rQ_\lambda,Q_{(\lambda,r)} \rangle=\langle Q_\lambda,Q_{(\lambda,0)}\rangle=2^{l(\lambda)}.
\end{equation}

The lemma then follows from a combination of equations \eqref{eqn:tr}, \eqref{eqn:trc1}, \eqref{eqn:trc2}, \eqref{eqn:trc3}, and the fact
$(-1)^{r-k}Q_{(\lambda,k,r-k)}=(-1)^{k}Q_{(\lambda,r-k,k)}$ for $0<k<r$
since $r$ is odd.
$\Box$

Note that if $\la$ is a strict partition, the partitions appeared on the right hand sides of
equations \eqref{eqn:dQ} and \eqref{eqn:multi pr} may not be strict any more. In fact, they are even not
necessarily partitions since negative parts might occur. In order to repeatedly using the above formula, it
is convenient to have the following

\begin{cor}\label{cor:pr-nota}
The derivative formula \eqref{eqn:dQ} and the multiplication formula \eqref{eqn:multi pr} hold for all partitions $\lambda$
which are not necessarily strict.
\end{cor}
This result might be known to experts. However we could not find a satisfactory reference for it. So we will give
a proof of this corollary in Appendix \ref{sec:prdr}.

In \cite{ASY}, Aokage, Shinkawa, and Yamada introduced  a different set of Virasoro operators $\{L_k'\}_{k\geq 1}$ and
computed their actions on Q-polynomials. In this paper we will only use the actions of the first two operators which
can be written as
\begin{align}
L_1':=\sum_{k \geq 1 \atop k \,\, odd} k t_k \frac{\partial}{\partial t_{k+2}}+\frac{1}{4}\frac{\partial^2}{\partial t_{1}\partial t_{1}},\hspace{20pt}
L_2' :=\sum_{k \geq 1 \atop k \,\, odd} k t_k \frac{\partial}{\partial t_{k+4}}+\frac{1}{2}\frac{\partial^2}{\partial t_{1}\partial t_{3}}.
            \label{eqn:L2'}
\end{align}
The coefficients of these operators are slightly different from the formulas given in \cite{ASY}. The reason for
this is that the definition of Q-polynomials in \cite{ASY} is slightly different from the definition we are using in
this paper. If we denote Q-polynomials in \cite{ASY} by $Q_\la^{\rm ASY}$, then
\[ Q_\la (\textbf{t})=(-1)^{\lceil l(\la)/2 \rceil} Q_\la^{\rm ASY} (2 \textbf{t}). \]
So the variables $t_k$ in \cite{ASY} is $2 t_k$ in this paper.

The following formulas were obtained in \cite{ASY} Theorem 2: For strict partitions $\la$,
\begin{align}\label{eqn:L'-formula}
L_1'Q_{\lambda}(\textbf{t})=\sum_{i=1}^{l(\lambda)} (\lambda_i-1) Q_{\lambda-2\epsilon_i}(\textbf{t}),\
L_2'Q_{\lambda}(\textbf{t})=\sum_{i=1}^{l(\lambda)} (\lambda_i-2) Q_{\lambda-4\epsilon_i}(\textbf{t}).
\end{align}
Since both sides of these equations are linear in Q-polynomials, the difference of signs for the definitions of
Q-polynomials does not affect these formulas as long as length of partitions appeared in both sides
of these equations do not change. In fact, the length of partition could change
only if $\la_i - 2$ or $\lambda_i-4$ becomes negative
for some $i$.
Since $\lambda$ is positive,  $\la_i - 2 < 0$ only if $\la_i=1$, so the corresponding coefficient in
the formula for $L_1'Q_{\lambda}$ becomes $0$. However, if $\la_i=3$ and $\la_j=1$ for some
$i$ and $j$,  $Q_{\lambda-4\epsilon_i}$ or $Q_{\lambda-4\epsilon_j}$ will be replaced by $Q_{\mu}$ for some partition $\mu$ with $l(\mu)=l(\la)-2$.
This will cause a change of sign
when comparing Q-polynomials in this paper and that in \cite{ASY}. However, the definitions for
strict partitions are also different in these two papers. Strict partitions in this paper have decreasing parts, while
strict partitions in \cite{ASY} have increasing parts. The order of the parts will determine whether
 $Q_{\lambda-4\epsilon_i}$ or $Q_{\lambda-4\epsilon_j}$ should contribute after removing negative parts,
 and their coefficients in $L_2'Q_{\lambda}$ happens
 to be either $1$ or $-1$. This will compensate the problem about the change of sign when lengths of partitions
 drop by $2$.  In a summary, equations \eqref{eqn:L'-formula} can be used in our setting without
 modifications.

\section{Computing $Q_{\lambda}(\delta_{k,3}/3)$}
\label{sec:dk3}

The definition of $\tau_{MM}$ as given by equation \eqref{eqn:tauMM}  depends on two special values of
Q-polynomials:
$Q_\lambda(\delta_{k,1})$ and $Q_\lambda(\delta_{k,3}/3)$.
It is well-known that
\begin{align}\label{eqn:hook}
Q_\lambda(\delta_{k,1})=\frac{2^{|\lambda|}}{\lambda!}\prod_{i<j}\frac{\lambda_i-\lambda_j}{\lambda_i+\lambda_j},
\end{align}
where $\la !:= \prod_{i=1}^{l(\la)} \la_{i} !$.
This formula can be interpreted as an analogue of the hook length formula (see, for example, equation (3.3) in \cite{Alex2020}).
In this section, we will give an explicit formula for $Q_{\lambda}(\delta_{k,3}/3)$ and study some basic properties
of these constants.

For any partition $\lambda=(\lambda_1, \cdots, \lambda_l)$, we define
\begin{equation}
 A_{\lambda} := Q_{\lambda}(\delta_{k,3}/3).
\end{equation}
For convenience, we may also allow $\la$ to have one negative part. In this case,
$A_{\la}$ should be interpreted as in Definition \ref{def:Q-}.
It follows from properties of $Q_{\la}$ that
$A_{\la}$ is skew symmetric with respect to permutations of $\la$ if $\la$ is weakly positive. Moreover
 $A_{(\lambda, 0)} = A_{\lambda}$ for all $\lambda$. We also have
 $A_{\la}=0$ if all parts of $\la$ are non-negative and $\la$ has two positive parts which are equal.

When dealing with $A_{\lambda}$, we often write $\lambda$ in the following standard form
\begin{equation} \label{eqn:Stform}
 \lambda=(3k_1, \cdots, 3k_p, 3m_1+1, \cdots, 3m_q+1, 3n_1+2, \cdots, 3n_r+2),
\end{equation}
where $k_i$, $m_i$, $n_i$ are non-negative integers.
For such $\lambda$, we also define
\begin{equation}
 \lambda_{[1]} := (3k_1, \cdots, 3k_p),  \hspace{10pt} \lambda_{[2]}:=(3m_1+1, \cdots, 3m_q+1, 3n_1+2, \cdots, 3n_r+2).
\end{equation}

\begin{thm} \label{thm:A}
For weakly positive $\lambda$ given by equation \eqref{eqn:Stform},
\[
A_{\lambda}
= \frac{\delta_{q,r} (-1)^{r(r-1)/2} 2^r (2/3)^{|\lambda|/3}}{\prod_{i=1}^p k_i ! \prod_{i=1}^q m_i ! \prod_{i=1}^r n_i !}
     \prod_{1 \leq i < j \leq p} \frac{k_i-k_j}{k_i+k_j} \,\, \frac{\prod_{1 \leq i < j \leq r} (m_i-m_j)(n_i-n_j)}{\prod_{i,j=1}^r (m_i+n_j+1)},
\]
where $\delta_{q,r}=1$ if $q=r$ and $\delta_{q,r}=0$ if $q \neq r$. In particular,
\begin{equation} \label{eqn:cong}
 A_{\la} = 0 {\rm \,\,\, if \,\,\,} q \neq r.
\end{equation}
\end{thm}
{\bf Proof}:
By equations \eqref{eqn:Ql1} and \eqref{eqn:Qwp}, the generating function for $A_{(k)}$
is given by
\[ \sum_{k=0}^{\infty} A_{(k)} \, z^k = \exp \left(\frac{2}{3} \, z^3 \right)
= \sum_{k=0}^{\infty} \left(\frac{2}{3}\right)^k \frac{z^{3k}}{k!}, \]
where $z$ is a formal parameter. Hence we have
\begin{eqnarray}
A_{(3k)} &=& \left( \frac{2}{3} \right)^{k} \frac{1}{k!}   \label{eqn:Ap1} \label{eqn:A3k}
\end{eqnarray}
and $A_{(k)}=0$ if $k/3$ is not an integer.
By equations  \eqref{eqn:Qwp} and \eqref{eqn:Ql2},
\[ A_{(m,n)} = A_{(m)} A_{(n)} + 2 \sum_{i=1}^n (-1)^i A_{(m+i)} A_{(n-i)}. \]
It follows that $A_{(m,n)}=0$ if $(m+n)/3$ is not an integer.
By using the elementary combinatorial identity
\begin{equation}
\sum_{i=a}^b (-1)^i \comb{n}{i} = (-1)^a \comb{n-1}{a-1} + (-1)^b \comb{n-1}{b} ,
\end{equation}
we obtain
\begin{eqnarray}
A_{(3k_1, 3k_2)} &=& \left( \frac{2}{3} \right)^{k_1+k_2} \frac{1}{k_1 ! k_2 !} \, \, \frac{k_1-k_2}{k_1+k_2} \label{eqn:Ap2}
\end{eqnarray}
for $(k_1, k_2) \neq (0,0)$, and
\begin{eqnarray}
A_{(3m+1, 3n+2)} &=&  \left( \frac{2}{3} \right)^{m+n+1} \frac{2}{m ! n !(m+n+1)}. \label{eqn:Al<3}
\end{eqnarray}
For any weakly positive partition $\la=(\la_1, \cdots, \la_l)$ with $l$ even, by equation \eqref{eqn:Qwp},
\[ A_{\la} = {\rm Pf}(A_{(\la_i, \la_j)})_{1 \leq i, j \leq l}. \]
In this formula, $A_{(0,0)}$ should be replaced by $0$ if $\la_i=0$ for some $i$.
It follows that $A_{\la}$ is alternating with respect to permutations of parts of $\la$ and
\begin{equation}
 A_{\lambda} = 0 {\rm \,\,\, if \,\,\,} |\lambda|/3 {\rm \,\,\, is \,\,\, not \,\,\, an \,\,\, integer}.
\end{equation}
For $\la$ given by equation \eqref{eqn:Stform}, the matrix $(A_{(\la_i, \la_j)})_{1 \leq i, j \leq l}$
has the following form
\[ \left( \begin{array}{ccc} B & 0 & 0 \\
                             0 & 0 & C \\
                             0 & -C^T & 0 \end{array} \right), \]
where
\[ B=(A_{(3k_i, 3k_j)})_{1 \leq i, j \leq p}, \hspace{20pt} C=(A_{(3 m_i+1, 3n_j+2)})_{1 \leq i \leq q, 1 \leq j \leq r}. \]
If $q \neq r$, the matrix $(A_{(\la_i, \la_j)})_{1 \leq i, j \leq l}$ is singular and $A_{\la}=0$.
If $q=r$,
\[ A_{\la} = {\rm Pf}(B) \cdot {\rm Pf} \left( \begin{array}{cc} 0 & C \\
                                                           -C^T & 0 \end{array} \right)
= {\rm Pf}(B) \cdot (-1)^{r(r-1)/2} \det(C). \]
By equation \eqref{eqn:Ap2},
\[ {\rm Pf}(B) = \frac{(2/3)^{\sum_{i=1}^p k_i}}{\prod_{i=1}^p k_i !} {\rm Pf} \left( \frac{k_i - k_j}{k_i+k_j} \right)_{1 \leq i,j \leq p}
= \frac{(2/3)^{\sum_{i=1}^p k_i}}{\prod_{i=1}^p k_i !} \prod_{1 \leq i < j \leq p} \frac{k_i - k_j}{k_i+k_j}, \]
where the last equality is Schur's Pfaffian identity (c.f. \cite[Proposition 2.2]{O}).

By equation \eqref{eqn:Al<3},
\begin{eqnarray}
 \det(C) &=& \frac{2^r (2/3)^{\sum_{i=1}^{r} (m_i+n_i+1)}}{\prod_{i=1}^r m_i! n_i!}
                 \det \left( \frac{1}{m_i+n_j+1} \right)_{1 \leq i, j \leq r} \nonumber \\
       &=& \frac{2^r (2/3)^{\sum_{i=1}^{r} (m_i+n_i+1)}}{\prod_{i=1}^r m_i! n_i!}
            \frac{\prod_{1 \leq i < j \leq r} (m_i-m_j)(n_i-n_j)}{\prod_{i,j=1}^r (m_i+n_j+1)},  \nonumber
\end{eqnarray}
where the last equality is Cauchy's identity (c.f. \cite[p38]{HJ}).
The theorem then follows by combining the above three equations.
$\Box$

\begin{rem} \label{rem:A}
By skew symmetry of $Q_{\la}$, we can replace each strict partition $\la$ in the definition of $\tau_{MM}$
by any permutation of $\la$. In particular, we can take $\la$ to be a partition with  positive distinct parts and has the standard form \eqref{eqn:Stform}. So Theorem \ref{thm:A}, together with equation \eqref{eqn:hook} (or even simpler equation \eqref{eqn:Q/Q}),
give a complete explicit description of coefficients in $\tau_{MM}$.
\end{rem}

\begin{cor} \label{cor:A1&2}
For $\lambda$ given by equation \eqref{eqn:Stform},
\[
A_{\lambda}= A_{\lambda_{[1]}} A_{\lambda_{[2]}}.
\]
\end{cor}

\begin{cor} \label{cor:A2}
For $\lambda$ given by equation \eqref{eqn:Stform},
\[
A_{2\lambda}=A_{\lambda} \cdot
\frac{(-1)^{r} (1/3)^{|\lambda|/3}}{\prod_{i=1}^p (2k_i-1)!! \prod_{i=1}^q (2m_i-1)!! \prod_{i=1}^r (2n_i+1)!!}.
\]
\end{cor}
{\bf Remark}: In \cite{MMNO}, a more general formula for the ratio $Q_{\la}/Q_{N \la}$ at special
points $t_{k}=\frac{r}{2} \delta_{k,r}$ and a factorization formula  were obtained using fermionic representation of Q-polynomials
(c.f. equations (7.6)  and (7.9) in \cite{MMNO}). Our method of using Pfaffian to prove Theorem \ref{thm:A} and Corollary \ref{cor:A2} is different from
that used in \cite{MMNO}. This method could be adapted to obtain more general results, which we omit here since they are not
needed in proving Theorem \ref{thm:MM=KW} and Theorem \ref{cor:KW as BKP}.

\begin{cor} \label{cor:MMhyper}
$\tau_{MM}(\textbf{t}/2)$ is a hypergeometric tau function of the BKP hierarchy.
\end{cor}
{\bf Proof}:
Since
\begin{align}\label{eqn:Q/Q}
\frac{Q_\lambda(\delta_{k,1})}{Q_{2\lambda}(\delta_{k,1})} = \prod_{j=1}^{l(\lambda)} (2 \lambda_j-1)!!
\end{align}
(c.f. \cite{Alex2020} Equation (4.16)), we can use Corollary \ref{cor:A2} to represent $\tau_{MM}$
in the following form:
\begin{equation} \label{eqn:Hyperg}
\tau_{MM}(\textbf{t}/2) = \sum_{\la \in DP} \theta_{\la} \, 2^{-l(\la)}  Q_{\la}(\textbf{t}/2) Q_{\la}(\textbf{t*}/2),
\end{equation}
where $\textbf{t*}=(t^*_1, t^*_3, t^*_5, \cdots)$ with $t^*_k = 2 \delta_{k,3}/3$ for all $k$, and for $\la$ given by a permutation of the partition in equation \eqref{eqn:Stform},
\[
\theta_{\la}= \left( \frac{\hbar}{48} \right)^{|\la|/3} (-1)^r
      \prod_{i=1}^p \frac{(6k_i-1)!!}{(2k_i-1)!!} \prod_{i=1}^q  \frac{(6m_i+1)!!}{(2m_i-1)!!} \prod_{i=1}^r \frac{(6n_i+3)!!}{(2n_i+1)!!} .
\]
For any non-negative integer $k$, define
\begin{eqnarray*}
\eta_{(3k)} &:=& \prod_{j=1}^k \frac{\hbar}{16} (6j-1)(6j-5), \\
\eta_{(3k+1)} &:=& c_1(\hbar) \prod_{j=1}^k \frac{\hbar}{16} (6j+1)(6j-1), \\
\eta_{(3k+2)} &:=& c_2(\hbar) \prod_{j=1}^k \frac{\hbar}{16} (6j+1)(6j-1)
\end{eqnarray*}
with any constants $c_1(\hbar)$ and $c_2(\hbar)$ such that
\[ c_1(\hbar) c_2(\hbar) = - \frac{\hbar}{16}. \]
For example, we can take $c_1(\hbar)=1$ and $c_2(\hbar) = - \frac{\hbar}{16}$, or
 $c_1(\hbar)=\left(- \frac{\hbar}{16}\right)^{1/3} $ and $c_2(\hbar) = \left(- \frac{\hbar}{16}\right)^{2/3}$.
Equivalently
\begin{eqnarray*}
\eta_{(3k)} &=& \left( \frac{\hbar}{48} \right)^{k} \frac{(6k-1)!!}{(2k-1)!!}, \\
\eta_{(3k+1)} &=& c_1(\hbar) \left( \frac{\hbar}{48} \right)^{k} \frac{(6k+1)!!}{(2k-1)!!}, \\
\eta_{(3k+2)} &=& c_2(\hbar) \frac{1}{3} \left( \frac{\hbar}{48} \right)^{k} \frac{(6k+3)!!}{(2k+1)!!}.
\end{eqnarray*}
Define
\[ \eta_{\la} := \prod_{i=1}^{l(\la)} \eta_{(\la_i)}. \]
Note that the only contributions to the right hand side of equation \eqref{eqn:Hyperg} are from $\la$ which are permutations
of partitions in the form of equation \eqref{eqn:Stform} with $q=r$. For such $\la$ we have
\[ \theta_{\la} = \eta_{\la}. \]
So equation \eqref{eqn:Hyperg} implies that $\tau_{MM}(\textbf{t}/2)$ is a hypergeometric tau function as defined in \cite{Or} (see also \cite{Alex2020}).
Note that if we take $c_1(\hbar) = \frac{\beta \hbar^{1/3}}{2}$ and $c_2(\hbar) = - \frac{ \hbar^{2/3}}{8 \beta}$, then
$\eta_{\la}$ is equal to $r_{\la}^{KW}$ used by Alexandrov in formulating his Conjecture 3 in \cite{Alex2020}.
$\Box$

\section{Virasoro constraints for $\tau_{MM}$}
\label{sec:VirMM}

In this section, we will prove $\tau_{MM}$ satisfies the {\it Virasoro constraints}
\[ L_k \tau_{MM} = 0\]
 for all $k \geq -1$, where $L_k$ is given by equation \eqref{eqn:Lk}.
This is equivalent to show
\[ \langle L_{k} \tau_{MM}, Q_{\mu} \rangle = 0 \]
for all strict partitions $\mu$.
Note that coefficients of $Q_{\la}$ in $\tau_{MM}$ are skew symmetric with respect to permutations
of $\la$. The following formula will be useful in proving Virasoro constraints for $\tau_{MM}$.

\begin{lem} \label{lem:inner-consis}
Let $f(\la)$ be any function which only depends on partitions $\la$ and not on $\textbf{t}$.
If $f(\la)$ is skew symmetric with respect to permutations of two parts of $\la$ which are not both $0$ and
$f((\la,0))=f(\la)$ for all partitions $\la$, then
 for any partition $\mu$,
\begin{equation} \label{eqn:inner-consis}
 \left\langle \sum_{\la \in DP} 2^{-l(\la)} f(\la) Q_{\la}, \,\, Q_{\mu} \right\rangle = f(\mu).
\end{equation}
\end{lem}
{\bf Proof}: Note that $f(\la) Q_{\la}$ is symmetric with respect to permutations of all parts of $\la$.
So in $\sum_{\la \in DP} 2^{-l(\la)} f(\la) Q_{\la}$ we can replace each $\la$ by any permutation of $\la$.
In particular, if $\mu$ is a permutation of some strict partition $\tilde{\mu}$, then
only $f(\tilde{\mu}) Q_{\tilde{\mu}}=f(\mu)Q_{\mu}$ term can contribute to the left hand side of equation \eqref{eqn:inner-consis}
by orthogonality of the Q-polynomials.
Hence we have
\[ \left\langle \sum_{\la \in DP} 2^{-l(\la)} f(\la) Q_{\la}, \,\, Q_{\mu} \right\rangle
= \langle 2^{-l(\mu)} f(\mu)Q_{\mu}, Q_{\mu} \rangle = f(\mu). \]

If $\mu$ has two equal positive parts, then $Q_\mu=0$ and $f(\mu)=0$ by skew symmetry of these two functions with respect to permutation of $\mu$. Hence both sides of equation \eqref{eqn:inner-consis} are zero. Therefore the lemma holds for
all positive partitions $\mu$.

We then prove the lemma by induction on the number of zero parts of $\mu$. Since both sides of
equation \eqref{eqn:inner-consis} are skew symmetric with respect to permutation of a zero part and
a positive part, if $\mu$ has a zero part, we permute $\mu$ to a partition of the form $(\mu', 0)$.
Since $Q_{(\mu',0)}=Q_{\mu'}$ and $f((\mu',0))=f(\mu')$, equation \eqref{eqn:inner-consis}
for $\mu$ is reduced to that of $\mu'$, which holds by induction hypothesis.
The lemma is thus proved.
$\Box$


Lemma \ref{lem:inner-consis} indicates that computing $\langle  \tau_{MM}, (L_{k})^\perp Q_{\mu} \rangle$
could be simpler than computing $\langle L_{k} \tau_{MM}, Q_{\mu} \rangle$ although these two quantities
are equal. Therefore instead of computing $L_{k} \tau_{MM}$, we will compute
$(L_{k})^\perp Q_{\mu} $ for all strict partitions $\mu$.
Since $L_{-1}$ and $L_2$ generate all $L_k$ for $k \geq -1$,
we only need to compute$(L_{k})^\perp Q_{\mu} $ for $k=-1, 2$.

First, we observe that operators  $L_{-1}$ and $L_2$ given by equation \eqref{eqn:VirOp}
and operators $L_1'$ and $L_2'$ given by equation \eqref{eqn:L2'} have the following relation:
\begin{lem} \label{lem:Lperp}
\begin{equation} \label{eqn:L-1perp}
(L_{-1})^\perp=\frac{1}{2}L_1'-\frac{1}{16}\frac{\partial^2}{\partial t_{1}\partial t_{1}}-\frac{t_1}{\hbar},
\end{equation}
\begin{equation} \label{eqn:L2perp}
(L_2)^\perp=\frac{1}{2}(L_2')^\perp+3t_1t_3-\frac{7t_7}{\hbar}.
\end{equation}
\end{lem}
{\bf Proof}:
Since \[ L_2= \frac{1}{2} L_2' + \frac{1}{4} \frac{\partial^2}{\partial t_1 \partial t_3}
              - \frac{1}{2 \hbar} \frac{\partial}{\partial t_7}, \]
This lemma follows from equation \eqref{eqn:ad-relation}.
$\Box$

Therefore, to compute $(L_{-1})^\perp \cdot Q_\lambda$ and $(L_2)^\perp \cdot Q_\lambda$,
  we need the following two lemmas.

\begin{lem} \label{lem:L2'adj} For any strict partition $\la=(\la_1, \cdots, \la_l)$,
\[(L_2')^\perp \cdot Q_\lambda=\sum_{i=1}^{l(\lambda)}(\lambda_i+2)Q_{\la+4\epsilon_i}+Q_{(\la,4)}-\frac{1}{2}Q_{(\la,3,1)}. \]
\end{lem}
{\bf Proof}:
By orthogonality of Q-polynomials,
\begin{equation} \label{eqn:L2'adjin}
(L_2')^\perp \cdot Q_\lambda
= \sum_{\mu \in DP} 2^{-l(\mu)} \, \langle (L_2')^\perp \cdot Q_\lambda, \,\, Q_{\mu} \rangle \,\, Q_{\mu}.
\end{equation}
By skew symmetry of Q-polynomials, $\langle (L_2')^\perp \cdot Q_\lambda, \, Q_{\mu} \rangle \, Q_{\mu}$ is symmetric with
respect to permutations of $\mu$. Hence we can replace each $\mu$ in the above equation by any permutation of $\mu$.
After permutation, $\mu$ may not be strict, but it must have distinct parts.

By equation \eqref{eqn:L'-formula}, for $\mu = (\mu_1, \cdots, \mu_{l(\mu)})$,
\begin{equation} \label{eqn:L2'ml}
 \langle (L_2')^\perp \cdot Q_\lambda, \,\, Q_{\mu} \rangle = \langle Q_\lambda, \,\, L_2' \cdot Q_{\mu} \rangle
=\sum_{j=1}^{l(\mu)} (\mu_j - 2) \langle Q_\lambda, \,\,   Q_{\mu - 4 \epsilon_j} \rangle.
\end{equation}
For this inner product to be non-zero, $\mu$ must have one of the following three forms.

{\bf Case (1)}, $\mu$ is a permutation of $\la + 4  \epsilon_i$ for some $i$ between $1$ and $l(\la)$. We may assume
$\mu = \la + 4  \epsilon_i$. By equation \eqref{eqn:L2'ml},
\[ \langle (L_2')^\perp \cdot Q_\lambda, \,\, Q_{\mu} \rangle
=\sum_{j=1}^{l(\la)} (\la_j + 4 \delta_{j,i} - 2) \langle Q_\lambda, \,\,   Q_{\la + 4  \epsilon_i - 4 \epsilon_j} \rangle.
\]
For $\langle Q_\lambda, \,\,   Q_{\la + 4  \epsilon_i - 4 \epsilon_j} \rangle \neq 0$,
$\la + 4  \epsilon_i - 4 \epsilon_j$ must be a permutation of $\la$. If $j \neq i$, this implies $\la_i + 4 = \la_j$ which is not possible
since $\mu = \la + 4  \epsilon_i$ must have distinct parts. Therefore we must have $j=i$ and
\begin{equation} \label{eqn:L2'c1}
 \langle (L_2')^\perp \cdot Q_\lambda, \,\, Q_{\mu} \rangle
=(\la_i + 2) \langle Q_\lambda, \,\,   Q_{\la} \rangle = 2^{l(\la)} (\la_i+2)= 2^{l(\mu)} (\la_i+2)
\end{equation}
 for $\mu = \la + 4  \epsilon_i$.

{\bf Case (2)}, $\mu$ is a permutation of $(\la, 4)$. We may assume
$\mu = (\la, 4)$. Since $\mu$ must have distinct parts,  $\la$ can not have parts equal to $4$.
By equation \eqref{eqn:L2'ml},
\[ \langle (L_2')^\perp \cdot Q_\lambda, \,\, Q_{\mu} \rangle
=\sum_{j=1}^{l(\la)} (\la_j  - 2) \langle Q_\lambda, \,\,   Q_{(\la-4\epsilon_j, 4)} \rangle
  + (4  - 2) \langle Q_\lambda, \,\,   Q_{(\la, 0)} \rangle.
\]
Since $\la$ does not have parts equal to $4$,
$\langle Q_\lambda, \,\,   Q_{(\la-4\epsilon_j, 4)} \rangle = 0$ for $1 \leq j \leq l(\la)$.
Moreover $Q_{(\la, 0)}=Q_{\la}$. So we have
\begin{equation} \label{eqn:L2'c2}
 \langle (L_2')^\perp \cdot Q_\lambda, \,\, Q_{\mu} \rangle
= 2 \langle Q_\lambda, \,\,   Q_{\la} \rangle = 2^{l(\la)+1}= 2^{l(\mu)}
\end{equation}
for $\mu = (\la, 4)$.

{\bf Case (3)}, $\mu$ is a permutation of $(\la, 3, 1)$. We may assume
$\mu = (\la, 3, 1)$. Since $\mu$ must have distinct parts,  $\la$ can not have
parts equal to $3$ or $1$.
By equation \eqref{eqn:L2'ml},
\[ \langle (L_2')^\perp \cdot Q_\lambda, \,\, Q_{\mu} \rangle
=\sum_{j=1}^{l(\la)} (\la_j  - 2) \langle Q_\lambda, \,\,   Q_{(\la-4\epsilon_j, 3,1)} \rangle
  + (3-2) \langle Q_\lambda, \,\,   Q_{(\la, -1,1)} \rangle.
\]
Since $\la$ does not have parts equal to $3$ and $1$,
$\langle Q_\lambda, \,\,   Q_{(\la-4\epsilon_j, 3,1)} \rangle=0$ for $1 \leq j \leq l(\la)$.
Moreover $Q_{(\la, -1,1)} = -2 Q_{\la}$. So we have
\begin{equation} \label{eqn:L2'c3}
 \langle (L_2')^\perp \cdot Q_\lambda, \,\, Q_{\mu} \rangle
= -2 \langle Q_\lambda, \,\,   Q_{\la} \rangle = -2^{l(\la)+1}= -2^{l(\mu)-1}
\end{equation}
for $\mu = (\la, 3, 1)$.

Combining equations \eqref{eqn:L2'adjin}, \eqref{eqn:L2'c1}, \eqref{eqn:L2'c2}, \eqref{eqn:L2'c3},
we obtain the desired formula. The lemma is thus proved.
$\Box$

{\allowdisplaybreaks
\begin{pro}\label{pro:DL}
For any strict partition $\lambda=(\la_1, \cdots, \la_l)$, we have
\begin{align}
(L_{-1})^\perp\cdot Q_{\lambda}
=&-\frac{1}{4}\sum_{i,j=1,\atop i\neq j}^{l(\lambda)} Q_{\lambda-\epsilon_i-\epsilon_j}
+\sum_{i=1}^{l(\lambda)}\frac{2\lambda_i-3}{4} Q_{\lambda-2\epsilon_i}  \nonumber \\
& -\frac{1}{\hbar}\bigg(\sum_{i=1}^{l(\lambda)} Q_{\lambda+\epsilon_i}
    + \frac{1}{2} Q_{(\lambda,1)}\bigg), \label{eqn:DL-1}\\
(L_{2})^\perp\cdot Q_{\lambda}
=&\sum_{i,j=1,\atop i\neq j}^{l(\lambda)} Q_{\la+3\epsilon_i+\epsilon_j}+\sum_{i=1}^{l(\lambda)}\frac{\lambda_i+4}{2} Q_{\la+4\epsilon_i}+\frac{1}{2}\sum_{i=1}^{l(\lambda)}Q_{(\la+3\epsilon_i,1)} \nonumber\\
&\ \ \ \ +\frac{1}{2}\sum_{i=1}^{l(\lambda)}Q_{(\la+\epsilon_i,3)}
-\frac{1}{2}\sum_{i=1}^{l(\lambda)}Q_{(\la+\epsilon_i,2,1)}
+Q_{(\la,4)}-\frac{1}{2}Q_{(\la,3,1)} \nonumber\\
&-\frac{1}{\hbar}\bigg(\sum_{i=1}^{l(\la)}Q_{\la+7\epsilon_i}
+\frac{1}{2} \sum_{r=0}^3 (-1)^r Q_{(\la,7-r, r)} \bigg). \label{eqn:DL2}
\end{align}
\end{pro}
}
{\bf Proof}:
Given a positive partition $\la$, by equation \eqref{eqn:dQ} and Corollary \ref{cor:pr-nota},
$ \frac{1}{2}\frac{\partial}{\partial t_{1}} Q_\la = \sum_{i=1}^{l(\lambda)}Q_{\lambda-\epsilon_i}, $
where each $\lambda-\epsilon_i$ on the right hand side is still a partition, i.e. all parts are non-negative.
By Corollary \ref{cor:pr-nota}, we can apply equation \eqref{eqn:dQ} again to take derivative with respect to $t_1$
on both sides of this equation and obtain 
\[  \frac{1}{4}\frac{\partial^2}{\partial t_{1}\partial t_{1}} Q_\lambda
        = \sum_{i,j=1}^{l(\lambda)} Q_{\lambda-\epsilon_i-\epsilon_j} . \]
By Corollary \ref{cor:pr-nota}, we can repeatedly using equations \eqref{eqn:multi pr} to multiply $Q_{\la}$ by $t_1$ and $t_3$
and  obtain a formula for computing $t_1 t_3 Q_{\lambda}$. 
Combining these results with Lemmas \ref{lem:Lperp}, \ref{lem:L2'adj}, and 
equations \eqref{eqn:L'-formula} and \eqref{eqn:multi pr}, we obtain the desired formulas.
$\Box$

\begin{rem}  \label{rem:DL}
It is straightforward to check that both sides of equations \eqref{eqn:DL-1} and \eqref{eqn:DL2}
do not change value if $\la$ is replaced by $(\la, 0)$. Hence these equations also hold for $\la=(\mu, 0)$
where $\mu$ is a strict partition.
\end{rem}

Now we are ready to compute $\langle L_{k} \, \tau_{MM}, \, Q_{\mu} \rangle$ for $k=-1, 2$.
They are essentially given by the following functions of $\mu$:
\begin{align}\label{eqn:Phi}
\Phi(\mu):=\sum_{i,j=1\atop i \neq j}^{l(\mu)}  \frac{A_{2\mu-2\epsilon_i-2\epsilon_j}}{(2\mu_i-1)(2\mu_j-1)}
- \sum_{i=1}^{l(\mu)} \frac{A_{2\mu - 4\epsilon_i}}{(2\mu_i-1)}
+ \frac{1}{4}\sum_{i=1}^{l(\mu)} (2\mu_i+1) A_{2\mu+2\epsilon_i}
+ \frac{1}{8} A_{(2\mu,2)},
\end{align}
and
{\allowdisplaybreaks
\begin{align}\label{eqn:Psi}
\Psi(\mu):=&\sum_{i=1}^{l(\mu)}\sum_{j=1, j\neq i}^{l(\mu)} (2\mu_j+1)(2\mu_i+5)_{[[3]]} A_{2\mu+6\epsilon_i+2\epsilon_j}
+\frac{1}{2}\sum_{i=1}^{l(\mu)}(\mu_i+4)(2\mu_i+7)_{[[4]]} A_{2\mu+8\epsilon_i} \nonumber\\
&+\frac{1}{2}\sum_{i=1}^{l(\mu)}(2\mu_i+5)_{[[3]]} A_{(2\mu+6\epsilon_i,2)}
+\frac{15}{2}\sum_{i=1}^{l(\mu)}(2\mu_i+1)A_{(2\mu+2\epsilon_i,6)} \nonumber\\
&-\frac{3}{2}\sum_{i=1}^{l(\mu)}(2\mu_i+1)A_{(2\mu+2\epsilon_i,4,2)}
+105A_{(2\mu,8)}
-\frac{15}{2}A_{(2\mu,6,2)} \nonumber\\
&-\frac{1}{16} \bigg( \sum_{i=1}^{l(\mu)}(2\mu_i+13)_{[[7]]} A_{2\mu+14\epsilon_i}
+\sum_{r=0}^3 \frac{(-1)^r (2r-1)!!(13-2r)!!}{2} A_{2(\mu,7-r,r)} \bigg),
\end{align}
}
where \[ x_{[[k]]}:= x(x-2)\cdots(x-2k+2) \]
 is the double falling factorial, and we set $(-1)!!=1$ for convenience.

\begin{thm}\label{thm:LtoPhi}
For all strict partitions $\mu$,
\begin{eqnarray*} 
\langle L_{-1} \, \tau_{MM}, \, Q_{\mu} \rangle
&=& -\frac{1}{4} \left(\frac{\hbar}{16}\right)^{(|\mu|-2)/3} \frac{Q_{\mu}(\delta_{k,1})}{Q_{2\mu}(\delta_{k,1})} \,\, \Phi(\mu), \\
\langle L_{2} \, \tau_{MM}, \, Q_{\mu} \rangle
&=&  \left(\frac{\hbar}{16}\right)^{(|\mu|+4)/3} \frac{Q_{\mu}(\delta_{k,1})}{Q_{2\mu}(\delta_{k,1})} \,\, \Psi(\mu).
\end{eqnarray*}
\end{thm}
{\bf Proof}:
For any partition $\la$, let
\[ a_{\la}:= \left(\frac{\hbar}{16}\right)^{|\la|/3} \frac{Q_{\la}(\delta_{k,1})}{Q_{2\la}(\delta_{k,1})}
\hspace{20pt} {\rm and \,\,\,\,\,\,\,\,\,}
 f_{MM}(\la) := \, a_{\la} A_{2 \la}.\]
Then $f_{MM}(\la)$ satisfies conditions for $f(\la)$ in Lemma \ref{lem:inner-consis} and
\[ \tau_{MM} = \sum_{\la \in DP} 2^{-l(\la)} f_{MM}(\la) \, Q_{\la}. \]
Since for $k=-1,2$,
\[ \langle L_{k} \, \tau_{MM}, \, Q_{\mu} \rangle = \langle \tau_{MM}, \, (L_{k})^\perp \, Q_{\mu} \rangle.\]
By Proposition \ref{pro:DL} and Lemma \ref{lem:inner-consis}, we have
\begin{eqnarray}
\langle L_{-1} \, \tau_{MM}, \, Q_{\mu} \rangle
&=& - \frac{1}{4} \sum_{i,j=1, \atop i\neq j}^{l(\mu)}  a_{\mu-\epsilon_i-\epsilon_j} A_{2(\mu-\epsilon_i-\epsilon_j)}
+\sum_{i=1}^{l(\mu)}\frac{2\mu_i-3}{4} a_{\mu-2\epsilon_i} A_{2(\mu-2\epsilon_i)}  \nonumber \\
&& -\frac{1}{\hbar} \sum_{i=1}^{l(\mu)} a_{\mu+\epsilon_i} A_{2(\mu+\epsilon_i)}
                            - \frac{1}{2\hbar} a_{(\mu,1)} A_{2(\mu,1)} ,  \label{eqn:L-1Qu}
\end{eqnarray}
and
{\allowdisplaybreaks
\begin{eqnarray}
\langle L_{2} \, \tau_{MM}, \, Q_{\mu} \rangle
&=&  \sum_{i,j=1, \atop i\neq j}^{l(\mu)} a_{\mu+3\epsilon_i+\epsilon_j} A_{2(\mu+3\epsilon_i+\epsilon_j)}
     +\sum_{i=1}^{l(\mu)}\frac{\mu_i+4}{2}  a_{\mu+4\epsilon_i} A_{2(\mu+4\epsilon_i)}  \nonumber \\
&&   +\frac{1}{2}\sum_{i=1}^{l(\mu)} a_{(\mu+3\epsilon_i,1)}  A_{2(\mu+3\epsilon_i,1)}
     +\frac{1}{2} \sum_{i=1}^{l(\mu)} a_{(\mu+\epsilon_i,3)} A_{2(\mu+\epsilon_i,3)}  \nonumber \\
&&  -\frac{1}{2} \sum_{i=1}^{l(\mu)} a_{(\mu+\epsilon_i,2,1)} A_{2(\mu+\epsilon_i,2,1)}
   + a_{(\mu,4)} A_{2(\mu,4)}
    -\frac{1}{2} a_{(\mu,3,1)} A_{2(\mu,3,1)} \nonumber \\
&& -\frac{1}{\hbar}\sum_{i=1}^{l(\mu)} a_{\mu+7\epsilon_i} A_{2(\mu+7\epsilon_i)}
                           - \sum_{r=0}^3 \frac{(-1)^r}{2\hbar} a_{(\mu,7-r,r)} A_{2(\mu,7-r,r)}.
                           \label{eqn:L2Qu}
\end{eqnarray}
}
 Note that since $\mu$ is a strict partition, $\mu-\epsilon_i-\epsilon_j$ with $i \neq j$ are still partitions, so there is
 no problem to use Lemma \ref{lem:inner-consis} to compute $\langle \tau_{MM}, Q_{\mu-\epsilon_i-\epsilon_j} \rangle$.
 However $\mu-2\epsilon_i$ might have a part equal to $-1$ if $\mu_{i}=1$ for some $i$. In this case
 we can not apply Lemma \ref{lem:inner-consis} to compute $\langle \tau_{MM}, Q_{\mu-2\epsilon_i} \rangle$ directly.
 Since $\mu$ is a strict partition, it has at most one part equal to $1$, Hence $\mu-2\epsilon_i$ can not have a
 pair of parts equal to $-1$ and $1$ respectively. So the formula
 \[ \langle \tau_{MM}, Q_{\mu-2\epsilon_i} \rangle = a_{\mu-2\epsilon_i} A_{2(\mu-2\epsilon_i)} \]
 still holds since both sides of this equation are $0$.

By equation \eqref{eqn:Q/Q},
\[ a_{\mu} = \left( \frac{\hbar}{16} \right)^{|\mu|/3} \prod_{j=1}^{l(\mu)} (2 \mu_j-1)!! \]
for all partition $\mu$.
The theorem is obtained by factorizing out $- \frac{a_{\mu}}{4}  \left( \frac{\hbar}{16} \right)^{-2/3}$
from right hand side of equation \eqref{eqn:L-1Qu} and $a_{\mu} \left( \frac{\hbar}{16} \right)^{4/3}$
from right hand side of equation \eqref{eqn:L2Qu}.
$\Box$

{\bf Proof of Theorem \ref{thm:MM=KW}}:
By Theorem \ref{thm:LtoPhi}, the $L_{-1}$ and $L_2$ constraints for $\tau_{MM}$ follow from two combinatorial
identities $\Phi(\mu)=0$ and $\Psi(\mu)=0$ for all strict partitions $\mu$. These identities can be proved
by induction on the length of $\mu$, and the details will be given in sections \ref{sec:Phi=0}
and \ref{sec:Psi=0} respectively.
 Since $L_{-1}$ and $L_2$ generate all $L_k$ with $k \geq -1$, this proves Theorem \ref{thm:MM=KW}.
$\Box$

{\bf Proof of Corollary \ref{cor:KW=MM}}:
By Theorem \ref{thm:A}, $Q_{\la}(\delta_{k,3}/3) = A_{\la}=0$ if $|\la|/3$ is not an integer.
Hence only non-negative integral powers of $\hbar$ could occur in $\tau_{MM}$.
Since $Q_{\la}$ is a homogeneous polynomial of degree $|\la|$,
the coefficient of $\hbar^m$ in $\tau_{MM}$ is a homogeneous polynomial of degree $3m$
for all $m \geq 0$.
This corresponds to the dimension constraint for $\tau_{KW}$ (see Section \ref{sec:KW}).
Moreover, at $\textbf{t}=0$, $\tau_{MM} = Q_{\emptyset}=1$. Hence
$\tau_{MM}$ and $\tau_{KW}$ have the same initial value.
By Theorem \ref{thm:MM=KW} and Kontsevich's theorem, both $\tau_{MM}$ and $\tau_{KW}$
satisfy the same Virasoro constraints.
Since the Virasoro constraints determine the tau functions up to a constant (cf. \cite{D} and \cite{VV}),
this completes the
proof of Corollary \ref{cor:KW=MM}.
$\Box$

\section{Combinatorial identities for constants $A_{\la}$}
\label{sec:Comb}

Recall $A_\la=Q_\la(\delta_{k,3}/3)$. The values of these constants are explicitly given in Theorem \ref{thm:A}.
In this section, we prove some combinatorial identities satisfied by these constants.
Since $\tau_{MM}$ is determined by these constants, properties of $A_\la$ will be important in the further study
of this tau function. 
In particular, the two identities $\Phi(\mu)=0$ and $\Psi(\mu)=0$, where $\Phi(\mu)$ and $\Psi(\mu)$ are 
  defined by equations \eqref{eqn:Phi} and \eqref{eqn:Psi}, will be proved in sections \ref{sec:Phi=0} and \ref{sec:Psi=0} respectively.
These identities were used in the proof of the $L_{-1}$ and $L_2$ constraints for $\tau_{MM}$  in section \ref{sec:VirMM}.

\subsection{Some elementary identities for $A_{\la}$}

In this subsection we collect
some useful identities for constants $A_\la$ which can be proved using basic properties of $Q_\la$.

$A_\lambda$ satisfies two recursion relations. The first recursion is
\begin{equation} \label{eqn:ARec1}
\frac{|\lambda|}{2} A_{\lambda} = \sum_{i=1}^{l(\lambda)} A_{\lambda-3 \ep_i}.
\end{equation}
This formula can be proved using equation \eqref{eqn:dQ}. Since we will not use this formula in this paper, we omit the proof.

The second recursion is
\begin{equation} \label{eqn:ARec2}
 A_{\lambda} = \sum_{i=2}^{l} (-1)^i A_{(\lambda_1, \lambda_i)} A_{\lambda^{\{1,i\}}}
\end{equation}
for any weakly positive partition $\lambda=(\lambda_1, \cdots, \lambda_l)$ with $l$ even.
This formula follows from similar properties for Pfaffian
(cf. Equation (2.4) in \cite{O}, see also Theorem 9.14 in \cite{HH}).
By skew symmetry of $A$, for any fixed $j$ we can use $\lambda_j$  to replace the role of $\lambda_1$ in formula \eqref{eqn:ARec2}.
The corresponding formula for weakly positive partition with even length is
\begin{equation} \label{eqn:ARec2r}
 A_{\lambda} = (-1)^{j-1}\sum_{i=1,\atop i\neq j}^{l} (-1)^{\tilde{i}(j)} A_{(\lambda_j, \lambda_i)} A_{\lambda^{\{j,i\}}},
\end{equation}
where
\begin{equation} \label{eqn:j-tilde}
\tilde{i}(j)= \left\{ \begin{array}{ll} i-1, & {\rm if \,\,\,} i<j, \\
                                       i, & {\rm if \,\,\,} i > j.
                                       \end{array} \right.
\end{equation}
We call the right hand side of equation \eqref{eqn:ARec2r} the expansion of $A_{\la}$ with respect
to the $j$-th part $\la_j$.

If $\la$ is a positive partition with odd length, we should replace $\lambda$ by $(\lambda,0)$ in equation~\eqref{eqn:ARec2r}.
In particular, we have
\begin{equation} \label{eqn:ARecOdd0}
 A_{\lambda} = \sum_{i=1}^{l} (-1)^{i+1} A_{(\lambda_i)} A_{\lambda^{\{i\}}}
\end{equation}
for positive partition with odd length. This formula is obtained by expanding $A_{(\la,0)}$ using
equation \eqref{eqn:ARec2r} with $j=l+1$.

In comparison to equation \eqref{eqn:ARecOdd0},  we have:
\begin{lem} \label{lem:even-1}
If $\la$ is a positive partition with even length,
\[ \sum_{i=1}^{l} (-1)^{i+1} A_{(\lambda_i)} A_{\lambda^{\{i\}}}=0. \]
\end{lem}
{\bf Proof}:
Since $\lambda^{\{i\}}$ is a positive partition of odd length, we can use equation \eqref{eqn:ARecOdd0} to expand
$A_{\lambda^{\{i\}}}$ for all $i=1, \cdots, l$. Then the left hand side of the above equation becomes
\[ \sum_{i=1}^{l} (-1)^{i+1} A_{(\lambda_i)} \sum_{j=1 \atop j \neq i}^l (-1)^{\tilde{j}(i)} A_{(\la_j)} A_{\lambda^{\{i, j\}}}
= \sum_{i,j=1 \atop i \neq j}^{l} (-1)^{i+\tilde{j}(i)+1} A_{(\lambda_i)} A_{(\la_j)} A_{\lambda^{\{i, j\}}} = 0,
\]
where the last equality follows from the fact that each summand is skew symmetric with respect to $i$ and $j$.
$\Box$

Occasionally we also need to use recursion for $A_{\la}$ where $\la$ is a partition with two parts equal to 0.
Hence we need the following
\begin{lem} \label{lem:odd0}
Let $\mu=(\mu_1, \cdots, \mu_l)$ be a weakly positive partition of odd length.
Then equation \eqref{eqn:ARec2r} holds for $\la=(\mu, 0)$ with $j=1, \cdots, l$.
It also holds for $\la$ with $j=l+1$ up to a sign.
\end{lem}
{\bf Proof}:
Since $A_{\la}$ is skew symmetric with respect to permutations of the first $l$ parts of $\la$,
to prove equation \eqref{eqn:ARec2r} holds for $\la=(\mu, 0)$ with $j=1, \cdots, l$,
we only need to prove it holds for $j=1$, i.e. equation \eqref{eqn:ARec2} holds.

If $\mu$ is positive, then $\la$ is weakly positive. So equation \eqref{eqn:ARec2} holds automatically
for $\la$. Therefore we only need to consider the case where $\mu$ has exactly one part equal to 0.

If $\mu_1=0$, the left hand side of equation \eqref{eqn:ARec2} is equal to $A_{\mu^{\{1\}}}$ since $l$ is odd.
The right hand side of equation \eqref{eqn:ARec2} is
\begin{eqnarray*}
\sum_{i=2}^l (-1)^i A_{(0, \mu_i)} A_{(\mu^{\{1, i\}},0)} + A_{(0,0)} A_{\mu^{\{1\}}}
&=& \sum_{i=2}^l (-1)^{i+1} A_{(\mu_i)} A_{\mu^{\{1, i\}}} +  A_{\mu^{\{1\}}} = A_{\mu^{\{1\}}},
\end{eqnarray*}
where the last equality follows from Lemma \ref{lem:even-1} applied to the positive partition
$\mu^{\{1\}}$. Hence equation \eqref{eqn:ARec2} holds for this case.

Similar, the right hand side of equation \eqref{eqn:ARec2r} with $j=l+1$ is $-A_{\mu^{\{1\}}}$. Hence it only holds up to a sign.

If $\mu_i=0$ for some $2 \leq i \leq l$ and all other parts of $\mu$ are positive, then $\mu^{\{i\}}$ is a positive partition with even length. By
equation \eqref{eqn:ARec2},
\begin{equation} \label{eqn:mu-i}
  A_{\mu^{\{i\}}} = \sum_{k=2 \atop k \neq i}^l (-1)^{\tilde{k}(i)+1} A_{(\mu_1, \mu_{k})} A_{\mu^{\{1,k,i\}}}.
\end{equation}
On the other hand, plugging $\la=(\mu,0)$ into the right hand side of equation \eqref{eqn:ARec2}, we obtain
\begin{equation} \label{eqn:mu+0r}
 \sum_{k=2 \atop k \neq i}^l (-1)^k  A_{(\mu_1, \mu_{k})} A_{(\mu^{\{1,k\}},0)} + (-1)^{i} A_{(\mu_1, 0)} A_{(\mu^{\{1,i\}},0)}
+ (-1)^{l+1} A_{(\mu_1, 0)} A_{\mu^{\{1\}}}.
\end{equation}
The last two terms in this expression cancel each other since
\[ A_{\mu^{\{1\}}} = (-1)^{l-i} A_{(\mu^{\{1, i\}},0)}. \]
Moreover,
\[ A_{(\mu^{\{1,k\}},0)} = (-1)^{l-\tilde{i}(k)} A_{\mu^{\{1,k, i\}}}.
\]
So expression \eqref{eqn:mu+0r} is equal to $(-1)^{l-i}$ times the right hand side of equation \eqref{eqn:mu-i}.
Since $A_{\la} = (-1)^{l-i} A_{\mu^{\{i\}}}$,
equation \eqref{eqn:ARec2} holds for $\la$. The lemma is thus proved.
$\Box$

Using recursion formula \eqref{eqn:ARec2r}, we can prove the following very useful formula:
\begin{lem} \label{lem:ni-1}
For weakly positive partition $\lambda$ given by equation \eqref{eqn:Stform},
\[\sum_{i=1}^r A_{\lambda-\ep_{i+p+q}}=0.\]
\end{lem}
{\bf Proof}:
If $r \neq q+2$, then $A_{\lambda-\ep_{i+p+q}}=0$ for all $1 \leq i \leq r$. So the lemma is trivial in this case.

Assume $r=q+2$. We prove this lemma by induction on $p$ and $q$.

When $p=q=0$, the lemma has the form
\begin{equation} \label{eqn:ASym2}
A_{(3n_1+1, \, 3n_2+2)} + A_{(3n_1+2, \, 3n_2+1)}=0
\end{equation}
for all non-negative integers $n_1$ and $n_2$,
which follows from equation \eqref{eqn:Al<3} and skew symmetry of $A$.

Let
\[ R(\lambda)= \sum_{i=1}^r A_{\lambda-\ep_{i+p+q}}. \]
By Corollary \ref{cor:A1&2},
\[ R(\la) = A_{\la_{[1]}} R(\la_{[2]}). \]
So we only need to consider the case $p=0$. In this case, $l:=p+q+r$ is even.

If $q >0$, we use equation \eqref{eqn:ARec2r} with $j=l$ to expand all $A_{\lambda-\ep_{i+p+q}}$ and obtain
\begin{eqnarray}
 R(\lambda) &=&  \sum_{i=1}^q (-1)^{p+i+1+l} A_{(3m_i+1, \, 3n_r+2)} R(\lambda^{\{p+i, l\}}) \nonumber \\
    &&  + \sum_{i=1}^{r-1} (-1)^{p+q+i+1+l} A_{\lambda^{\{p+q+i, l\}}} \cdot \{ A_{(3n_i+1, \, 3n_r+2)} + A_{(3n_i+2, \, 3n_r+1)} \}.
    \label{eqn:ni-1Rec}
\end{eqnarray}
 Note in this calculation, we have used the fact
\[ A_{(\la_j, \, 3n_r+2)} = 0 \]
unless $\la_j = 3 a+1$ for some integer $a$.
The second term in the right hand side of equation~\eqref{eqn:ni-1Rec} is $0$ due to equation \eqref{eqn:ASym2}.
Since the number of parts in $\lambda^{\{p+i, l\}}$ which are equal to $1 {\rm \,\, mod}(3)$ is equal to $q-1$, by induction on $q$,
we only need to prove $R(\lambda)=0$ for $q=0$.
Hence the lemma follows from equation \eqref{eqn:ASym2}.
$\Box$

Setting $t_{k}=\delta_{k,3}/3$ for all $k$ and $r=1$ or $r=3$ in equation \eqref{eqn:multi pr}, we obtain the following
useful formulas:
\begin{equation} \label{eqn:Ala1}
A_{(\la, 1)} = - 2 \sum_{i=1}^{l(\la)} A_{\la+\ep_i}
\end{equation}
and
\begin{equation} \label{eqn:Ala21}
A_{(\la,2,1)}=2\sum_{i=1}^{l(\la)}A_{\la+3\epsilon_i}+A_{(\la,3)}-2A_\la
\end{equation}
for all $\la$. As a corollary of equation \eqref{eqn:Ala1}, we obtain the following lemma.

\begin{lem}\label{lem:qi+1}
	For $\lambda$ given by equation $(\ref{eqn:Stform})$,
	\[\sum_{i=1}^qA_{\lambda+\epsilon_{i+p}}=0.\]
\end{lem}
{\bf Proof}:
By equation \eqref{eqn:cong}, if $q\neq r+2$, then $A_{\lambda+\epsilon_{j+p}}=0$ for all $1\leq j\leq q$. So we can assume $q=r+2$.
By equation \eqref{eqn:Ala1},
\[ - \frac{1}{2} \, A_{(\la, 1)} =  \sum_{j=1}^p A_{\lambda+\epsilon_{j}} + \sum_{i=1}^q A_{\lambda+\epsilon_{i+p}}
              +\sum_{k=1}^r A_{\lambda+\epsilon_{k+p+q}}. \]
Since $q=r+2$, except terms in the middle summation on the right hand side, all other terms in above
equation are zero by equation \eqref{eqn:cong}. This proves the lemma.
$\Box$

\subsection{The identity $\Phi(\mu)=0$}
\label{sec:Phi=0}

In this subsection, we prove the identity $\Phi(\mu)=0$ which was used in the proof of the  $L_{-1}$-constraint for $\tau_{MM}$
in section \ref{sec:VirMM}. During the proof, we also obtain other non-trivial identities for $A_\la$ in
Lemma \ref{lem:Phi2'} and Propositions \ref{prop:Phi2} and \ref{prop:Phi3}. The main result of this subsection 
is the following
\begin{thm} \label{thm:MML-1}
Let $\Phi(\mu)$ be the function of $\mu$ defined by equation \eqref{eqn:Phi}. Then $\Phi(\mu)=0$ for
 all strict partitions $\mu$.
\end{thm}
{\bf Proof}:
Note that if $\mu$ is a strict partition, then for all $A_{\la}$ appeared in
the right hand side of equation \eqref{eqn:Phi}, $\la$ can not have two components which add up to $0$.
Therefore $\Phi(\mu)$ is skew symmetric with respect to permutations of $\mu$.
So to prove $\Phi(\mu)=0$ for a strict partition $\mu$, it is equivalent to prove
$\Phi(\nu)=0$  where $\nu$ is a permutation of $\mu$.

Without loss of generality, we may assume $\mu$ has the standard form
\begin{equation} \label{eqn:muStform}
 \mu=(3k_1, \cdots, 3k_p, 3m_1+1, \cdots, 3m_q+1, 3n_1+2, \cdots, 3n_r+2),
\end{equation}
where $k_i$, $m_i$, $n_i$ are non-negative integers.
Note that all $k_i$ must be positive if $\mu$ is a permutation of a strict partition, and at most one $k_i$ can be $0$
if $\mu$ is weakly positive.

In the rest part of this proof,
we will assume $\mu$ is a permutation of a strict partition.

We can use Corollary \ref{cor:A2} to represent $A_{2\la}$ in terms of $A_{\la}$ for all partitions $\la$ appeared in
$\Phi(\mu)$. Then after divided by constant
\[  \frac{(-1)^{r} (1/3)^{(|\mu|-2)/3}}{\prod_{i=1}^p (2k_i-1)!! \prod_{i=1}^q (2m_i-1)!! \prod_{i=1}^r (2n_i+1)!!}, \]
equation $\Phi(\mu)=0$  becomes
\[ \delta_{r, q+4} \cdot \frac{1}{9} \, \Phi_1 (\mu) + \delta_{r, q-2} \cdot \Phi_2(\mu) + \delta_{r, q+1} \cdot \Phi_3(\mu) = 0, \]
where $\Phi_i(\mu)$ are defined by
\begin{equation}
	\Phi_1(\mu) := \sum_{i,j=1 \atop i\neq j}^r A_{\mu-\epsilon_{i+p+q}-\epsilon_{j+p+q}},
\end{equation}
\begin{eqnarray}
\Phi_2(\mu)
&:=&  \sum_{i,j=1 \atop i\neq j}^p\frac{A_{\mu-\epsilon_{i}-\epsilon_{j}}}{(6k_i-1)(6k_j-1)}
      +\sum_{i,j=1 \atop i\neq j}^q\frac{A_{\mu-\epsilon_{i+p}-\epsilon_{j+p}}}{(6m_i+1)(6m_j+1)}  \nonumber \\
&& -2\sum_{i=1}^p\sum_{j=1}^q\frac{A_{\mu-\epsilon_{i}-\epsilon_{j+p}}}{(6k_i-1)(6m_j+1)}
   +\sum_{i=1}^q\frac{A_{\mu-2\epsilon_{i+p}}}{6m_i+1}
   -\frac{1}{4}\sum_{i=1}^q A_{\mu+\epsilon_{i+p}},  \label{eqn:Phi2}
\end{eqnarray}
and
\begin{eqnarray}
\Phi_3(\mu)
&:=&  \frac{2}{3} \sum_{i=1}^p\sum_{j=1}^r\frac{A_{\mu-\epsilon_{i}-\epsilon_{j+p+q}}}{6k_i-1}
     -\frac{2}{3} \sum_{i=1}^q\sum_{j=1}^r\frac{A_{\mu-\epsilon_{i+p}-\epsilon_{j+p+q}}}{6m_i+1} \nonumber \\
&& -\sum_{i=1}^p\frac{2k_i-1}{6k_i-1} A_{\mu-2\epsilon_{i}}
   +\frac{1}{3}\sum_{i=1}^r A_{\mu-2\epsilon_{i+p+q}}
   +\frac{1}{12}\sum_{i=1}^p(6k_i+1) A_{\mu+\epsilon_{i}}  \nonumber \\
&& -\frac{1}{12}\sum_{i=1}^r(6n_i+5) A_{\mu+\epsilon_{i+p+q}}
   +\frac{1}{24} A_{(\mu,1)}. \label{eqn:Phi3}
\end{eqnarray}

Applying Lemma \ref{lem:ni-1} for $\la=\mu-\epsilon_{j+p+q}$ and then summing over all $j=1, \cdots, r$, we have
\[ \Phi_1(\mu)=0\]
for all $\mu$.

We will prove $\Phi_2(\mu)=0$ and $\Phi_3(\mu)=0$ in Propositions \ref{prop:Phi2} and \ref{prop:Phi3} respectively.
This will complete the proof of Theorem \ref{thm:MML-1}.
$\Box$

To prove $\Phi_2(\mu)=0$, we need the following
\begin{lem} \label{lem:Phi2'}
For $\la$ given by equation \eqref{eqn:Stform}, define
\begin{eqnarray*}
\Phi_2'(\la) &:=&  \sum_{i=1}^q \frac{(-1)^i (2/3)^{m_i}}{m_i ! (6m_i+1)}
                   \bigg\{ \sum_{j=1 \atop j \neq i}^q \frac{ A_{(\la-\ep_{p+j})^{\{p+i\}}}}{6m_j+1}
                           - \sum_{j=1}^p \frac{ A_{(\la-\ep_j)^{\{p+i\}}}}{6k_j-1} \bigg\}.
\end{eqnarray*}
We have
\[ \Phi_2'(\la)=0 \]
for all $\la$ which is a permutation of a strict partition.
\end{lem}
{\bf Proof}:
If $q \neq r+2$, then every term in $\Phi_2'(\la)$ is equal to 0 by equation \eqref{eqn:cong}.
So we may assume $q = r+2$. We prove this lemma by induction on $p$ and $r$.

If $p=r=0$, then $q=2$ and
\begin{equation} \label{eqn:Phi2'base}
 \Phi_2'(\la)= \frac{1}{(6m_1+1)(6m_2+1)} \cdot \left\{ - \frac{(2/3)^{m_1} A_{(3m_2)}}{m_1 !} +
         \frac{(2/3)^{m_2} A_{(3m_1)}}{m_2 !} \right\} = 0,
\end{equation}
where the last equality follows from equation \eqref{eqn:Ap1}.

If $p$ is a positive odd integer, then partitions appeared in the right hand side of the definition of $\Phi_2'(\la)$
are weakly positive of even length.
We can apply equation \eqref{eqn:ARec2} to each term in $\Phi_2'(\la)$ and obtain
\begin{eqnarray}
\Phi_2'(\la) &=& \sum_{i=1}^q \frac{(-1)^i (2/3)^{m_i}}{m_i ! (6m_i+1)}
                       \sum_{j=1 \atop j \neq i}^q (-1)^{p+\tilde{j}(i)+1} A_{\la^{\{1, p+i, p+j\}}}
                           \left\{ \frac{A_{(3k_1, 3 m_j)}}{6m_j+1} - \frac{A_{(3k_1-1, 3 m_j+1)}}{6k_1-1} \right\}
                           \nonumber \\
              &&  + \sum_{i=2}^p (-1)^i A_{(3k_1, 3k_i)} \Phi_2'(\la^{\{1, i\}}), \label{eqn:Phi2'Rec}
\end{eqnarray}
where $\tilde{j}(i)$ is defined in equation \eqref{eqn:j-tilde}.

By equations \eqref{eqn:Ap2} and \eqref{eqn:Al<3},
\begin{equation} \label{eqn:A3k3m}
\frac{A_{(3k, 3 m)}}{6m+1} - \frac{A_{(3k-1, 3 m+1)}}{6k-1} = \frac{(2/3)^{k+m} (6k+1)}{k! m! (6m+1)(6k-1)}
\end{equation}
for any integers $k \geq 1$ and $m \geq 0$. Hence the first term on the right hand side of equation~\eqref{eqn:Phi2'Rec}
is equal to
\[ (-1)^{p+1}  \frac{ (2/3)^{k_1} (6k_1+1)}{k_1 !(6k_1-1)}
   \sum_{i,j=1 \atop j \neq i}^q  (-1)^{i+\tilde{j}(i)} \frac{ (2/3)^{m_i+m_j} A_{\la^{\{1, p+i, p+j\}}}}{m_i ! m_j ! (6m_i+1) (6m_j+1)},
\]
which is equal to $0$ since each term in the summation is skew symmetric with respect to $i$ and $j$.
Hence we have
\begin{eqnarray}
\Phi_2'(\la) &=&  \sum_{i=2}^p (-1)^i A_{(3k_1, 3k_i)} \Phi_2'(\la^{\{1, i\}}). \label{eqn:Phi2'Rec2}
\end{eqnarray}

When $p$ is a positive even integer, we need to replace each partition $\nu$ appeared in the right hand side of
the definition of $\Phi_2'(\la)$ by $(\nu, 0)$ and use Lemma \ref{lem:odd0} when applying equation \eqref{eqn:ARec2}.
All calculations are similar to the odd $p$ case except that an extra term $A_{(3k_1,0)} \Phi_2'(\la^{\{1\}})$ should be added
to the right hand sides of equations \eqref{eqn:Phi2'Rec} and \eqref{eqn:Phi2'Rec2}.

Since the numbers of parts in $\la^{\{1\}}$ and $\la^{\{1, i\}}$ which could be divided by $3$ are less than $p$,
by induction, we only need to prove $\Phi_2'(\la)=0$ for $p=0$.

Let $l:=p+q+r$.
If $r>0$ and $p=0$, we first replace each partition $\nu=(\nu_1, \cdots, \nu_{l-1})$ appeared in the right hand side of
the definition of $\Phi_2'(\la)$ by $(\nu, 0)$, then
apply equation \eqref{eqn:ARec2r} with $j=l-1$  to each term in $\Phi_2'(\la)$.
Note that $\nu_{l-1}=3n_r+2$ in each term and $A_{(3 n_r+2,0)}=0$. Hence the effect of adding $0$ to $\nu$ could be ignored in the final expansion formula. After the expansion,
we obtain
\begin{eqnarray}
\Phi_2'(\la)= &=& \sum_{i=1}^q (-1)^{i} A_{(3m_i+1, 3 n_r+2)} \Phi_2'(\la^{\{i, l\}}).
                            \nonumber
\end{eqnarray}
Since the number of parts in $\la^{\{i, l\}}$ which are equal to $2 \,\,{\rm mod}(3)$ is $r-1$,
we only need to prove $\Phi_2'(\la)=0$ for the case $r=p=0$, which we have considered in equation \eqref{eqn:Phi2'base}.
The proof of the lemma is thus finished.
$\Box$

We are now ready to prove
\begin{pro} \label{prop:Phi2}
For $\Phi_2(\mu)$ defined by equation \eqref{eqn:Phi2}, we have
\[ \Phi_2(\mu)=0 \]
for all $\mu$ which is a permutation of a strict partition.
\end{pro}
{\bf Proof}:
Assume $\mu$ is given by equation \eqref{eqn:muStform}.
If $q \neq r+2$, then every term in $\Phi_2(\mu)$ is equal to 0 by equation \eqref{eqn:cong}.
So we may assume $q = r+2$. In this case the length of $\mu$ is even if and only if $p$ is even.

Note that the right hand side of equation \eqref{eqn:Phi2} could have
a term containing $A_{\nu}$ with some parts of $\nu$ equal to $-1$. This could occur when
$\mu$ has a part equal to $1$. However, since $\mu$ is a permutation of a strict partition,
$\nu$ can not have simultaneously one part equal to $-1$ while another part equal to $1$.
So the corresponding $A_{\nu}=0$.
Hence the presence of a possibly negative part in $\nu$ does not affect the calculations.

We prove this proposition by induction on $p$ and $r$. If $p=r=0$, then $q=2$ and
\begin{eqnarray} \label{eqn:Phi2base}
 \Phi_2(\mu)&=& \Phi_2((3m_1+1, 3m_2+1)) \nonumber \\
 &=&  \frac{2 A_{(3m_1, 3m_2)}}{(6m_1+1)(6m_2+1)}  + \frac{A_{(3m_1-1, \, 3m_2+1)}}{6m_1 +1} +
         \frac{A_{(3m_1+1, \, 3m_2-1)}}{6m_2 +1} \nonumber \\
 && - \frac{A_{(3m_1+2, \, 3m_2+1)}}{4} - \frac{A_{(3m_1+1, \, 3m_2+2)}}{4} \nonumber \\
 &=& 0 ,
\end{eqnarray}
where the last equality is obtained by straightforward calculation using equations \eqref{eqn:Ap2} and \eqref{eqn:Al<3}.

If $p$ is a positive even integer,
we apply equation \eqref{eqn:ARec2} to each term in $\Phi_2(\mu)$ and obtain
\begin{eqnarray}
\Phi_2(\mu) &=& 2 \sum_{i=2}^p \sum_{j=1}^q  \frac{(-1)^{p+j}}{6k_i-1}
                         A_{(\mu-\ep_i)^{\{1, p+j\}}}
                           \left\{- \frac{A_{(3k_1, 3 m_j)}}{6m_j+1} + \frac{A_{(3k_1-1, 3 m_j+1)}}{6k_1-1} \right\}
                           \nonumber \\
              && + 2 \sum_{i,j=1 \atop i \neq j}^q   \frac{(-1)^{p+i}}{6m_j+1}
                         A_{(\mu-\ep_{p+j})^{\{1, p+i\}}}
                           \left\{ \frac{A_{(3k_1, 3 m_i)}}{6m_i+1} - \frac{A_{(3k_1-1, 3 m_i+1)}}{6k_1-1} \right\}
                           \nonumber \\
              &&  + \sum_{i=2}^p (-1)^i A_{(3k_1, 3k_i)} \Phi_2(\mu^{\{1, i\}}). \label{eqn:Phi2RecE}
\end{eqnarray}
After using equation \eqref{eqn:A3k3m} to simplify the first two terms in the above equation, we have
\begin{eqnarray}
\Phi_2(\mu) &=& 2 (-1)^p \frac{(2/3)^{k_1} (6k_1+1)}{ k_1 ! (6k_1-1)} \Phi_2'(\mu^{\{1\}})
              + \sum_{i=2}^p (-1)^i A_{(3k_1, 3k_i)} \Phi_2(\mu^{\{1, i\}}), \label{eqn:Phi2RecE2}
\end{eqnarray}
where $\Phi_2'(\mu^{\{1\}})=0$ by Lemma \ref{lem:Phi2'}. Hence we have
\begin{eqnarray}
\Phi_2(\mu) &=& \sum_{i=2}^p (-1)^i A_{(3k_1, 3k_i)} \Phi_2(\mu^{\{1, i\}}) \label{eqn:Phi2RecE3}.
\end{eqnarray}

If $p$ is a positive odd integer, we need first replace each partition $\nu$ appeared in
the right hand side of equation \eqref{eqn:Phi2} by $(\nu, 0)$. Then
we can apply equation \eqref{eqn:ARec2} to expand each term in $\Phi_2(\mu)$. The calculations
are similar to the $p$ even case except that an extra term
$A_{(3k_1, 0)} \Phi_2 (\mu^{\{1\}})$ should be added to equations
\eqref{eqn:Phi2RecE}, \eqref{eqn:Phi2RecE2}, and \eqref{eqn:Phi2RecE3}.

Since the numbers of parts in $\mu^{\{1\}}$ and $\mu^{\{1, i\}}$ which can be divided by 3 are less than $p$,
the proposition is reduced to the case of $p=0$.

If $r>0$ and $p=0$, applying equation \eqref{eqn:ARec2r} with $j=l:=q+r$ to each term in $\Phi_2(\mu)$, we obtain
\begin{eqnarray}
\Phi_2(\mu) &=& \sum_{i=1}^q (-1)^{i+1} A_{(3m_i+1, 3 n_r+2)} \Phi_2(\mu^{\{i, l\}}).
                            \nonumber
\end{eqnarray}
Since the number of parts in $\mu^{\{i, l\}}$ which are equal to $2 \,\,{\rm mod}(3)$ is $r-1$,
this equation reduces the proof to the case $r=p=0$, which we have considered in equation \eqref{eqn:Phi2base}.
The proof of the proposition is thus finished.
$\Box$

\begin{pro} \label{prop:Phi3}
For $\Phi_3(\mu)$ defined by equation \eqref{eqn:Phi3}, we have
\[ \Phi_3(\mu)=0 \]
for all $\mu$ which is a permutation of a strict partition.
\end{pro}
{\bf Proof}:
Assume $\mu$ is given by equation \eqref{eqn:muStform}.
If $r \neq q+1$, then every term in $\Phi_3(\mu)$ is equal to 0 by equation \eqref{eqn:cong}.
So we may assume $r = q+1$. In this case, the length of $\mu$ is even if and only if $p$ is odd.

By equation \eqref{eqn:Ala1}, the last term in $\Phi_3(\mu)$ is
\[ \frac{1}{24} A_{(\mu,1)} = - \frac{1}{12} \sum_{i=1}^{p+q+r} A_{\mu+\ep_i}.
\]
By equation \eqref{eqn:cong},  $A_{\mu+\ep_{p+j}}=0$ for all $1\leq j \leq q$ since $r = q+1$. Hence we can
get rid of the term $A_{(\mu,1)}$ in $\Phi_3(\mu)$ and obtain
\begin{eqnarray}
\Phi_3(\mu)
&=&  \frac{2}{3} \sum_{i=1}^p\sum_{j=1}^r\frac{A_{\mu-\epsilon_{i}-\epsilon_{j+p+q}}}{6k_i-1}
     -\frac{2}{3} \sum_{i=1}^q\sum_{j=1}^r\frac{A_{\mu-\epsilon_{i+p}-\epsilon_{j+p+q}}}{6m_i+1} \nonumber \\
&& -\sum_{i=1}^p\frac{2k_i-1}{6k_i-1} A_{\mu-2\epsilon_{i}}
   +\frac{1}{3}\sum_{i=1}^r A_{\mu-2\epsilon_{i+p+q}}
   +\frac{1}{2}\sum_{i=1}^p k_i A_{\mu+\epsilon_{i}}  \nonumber \\
&& -\frac{1}{2}\sum_{i=1}^r(n_i+1) A_{\mu+\epsilon_{i+p+q}}. \label{eqn:Phi3v2}
\end{eqnarray}
Since $\mu$ is a permutation of a strict partition, all partitions occurred in the right hand side of
this equation are weakly positive.

We now prove this proposition by induction on $p$ and $q$.
If $p=q=0$, then $r=1$ and
\begin{equation} \label{eqn:Phi3basep0}
 \Phi_3(\mu)= \Phi((3n+2))=\frac{1}{3}A_{(3n)}-\frac{n+1}{2} A_{(3n+3)} = 0,
\end{equation}
where the last equality follows from equation \eqref{eqn:Ap1}.

If $p=1, q=0$, then $r=1$ and
\begin{eqnarray}
\Phi_3(\mu)&=& \Phi_3((3k, 3n+2)) \nonumber \\
&=& \frac{2}{3} \frac{A_{(3k-1, 3n+1)}}{6k-1} - \frac{2k-1}{6k-1} A_{(3k-2, 3n+2)} \nonumber \\
&&    + \frac{1}{3} A_{(3k, 3n)} + \frac{k}{2} A_{(3k+1, 3n+2)} - \frac{n+1}{2} A_{(3k, 3n+3)} \nonumber \\
&=& 0, \label{eqn:Phi3base}
\end{eqnarray}
where the last equality follows from straightforward calculations  using  equations \eqref{eqn:Ap2} and \eqref{eqn:Al<3}.

If $p$ is a positive odd integer, we apply equation \eqref{eqn:ARec2} to expand each term in $\Phi_3(\mu)$ and obtain
\begin{eqnarray}
\Phi_3(\mu) &=& \frac{2}{3} \sum_{i=1}^q  (-1)^{p+i}
                           \left\{ \frac{A_{(3k_1-1, 3 m_i+1)}}{6k_1-1} - \frac{A_{(3k_1, 3 m_i)}}{6m_i+1} \right\}
                           \sum_{j=1}^r  A_{(\mu-\ep_{p+q+j})^{\{1, p+i\}}}
                           \nonumber \\
              && +  \sum_{i=1}^r   (-1)^{p+q+i} A_{\mu^{\{1, p+q+i\}}} \Phi_3((3k_1, 3n_i+2))
                           \nonumber \\
              &&  + \sum_{i=2}^p (-1)^i A_{(3k_1, 3k_i)} \Phi_3(\mu^{\{1, i\}}). \label{eqn:Phi3Rec}
\end{eqnarray}
The first term on the right hand side of the above equation is equal to $0$
after applying Lemma \ref{lem:ni-1} to $\la=\mu^{\{1, p+i\}}$ for each $i=1, \cdots, q$.
The second term is also equal to $0$ by equation
\eqref{eqn:Phi3base}. Hence we have
\begin{eqnarray}
\Phi_3(\mu) &=&  \sum_{i=2}^p (-1)^i A_{(3k_1, 3k_i)} \Phi_3(\mu^{\{1, i\}}). \label{eqn:Phi3Rec2}
\end{eqnarray}

If $p$ is a positive even integer, we first replace each partition $\nu$ occurred in the right hand side
of equation \eqref{eqn:Phi3v2} by $(\nu,0)$, then apply equation \eqref{eqn:ARec2} to expand each term in $\Phi_3(\mu)$.
The calculations are similar to $p$ odd case except that an extra term $A_{(3k_1,0)} \Phi_3(\mu^{\{1\}})$
should be added to the right hand sides of equations \eqref{eqn:Phi3Rec} and \eqref{eqn:Phi3Rec2}.
By induction, the proposition is then reduced to the case of $p=0$.

If $q > 0$ and $p=0$, we first replace each partition $\nu$ occurred in the right hand side
of equation \eqref{eqn:Phi3v2} by $(\nu,0)$, then apply equation \eqref{eqn:ARec2r} for $j=l:=q+r$
to expand each term in $\Phi_3(\mu)$. We obtain
{\allowdisplaybreaks
\begin{eqnarray}
\Phi_3(\mu)
&=&  - \frac{2}{3} \sum_{i=1}^q \sum_{j=1}^{r-1}  \frac{(-1)^{q+j}}{6m_i+1} A_{(\mu-\ep_{i})^{\{q+j, l\}}}
              \left\{ A_{(3n_j+1, 3n_r+2)} + A_{(3n_j+2, 3n_r+1)} \right\}   \nonumber  \\
&&   + A_{\mu^{\{l\}}} \Phi_3((3n_r+2))
     + \sum_{i=1}^q (-1)^{i} A_{(3m_i+1, 3 n_r+2)} \Phi_3(\mu^{\{i, l\}}).
\end{eqnarray}
}
 The first term on the right hand side of the above equation vanish due to
equation \eqref{eqn:ASym2}, the second term vanishes due to equation \eqref{eqn:Phi3basep0}. By induction,
the proposition is reduced to the case $q=p=0$, which has been considered in equation \eqref{eqn:Phi3basep0}.

The proof of the proposition is thus finished. This also completes the proof of Theorem \ref{thm:MML-1}.
$\Box$

\subsection{The identity $\Psi(\mu)=0$}
\label{sec:Psi=0}

In this subsection we prove the identity $\Psi(\mu)=0$ which was used in the proof of the $L_2$-constraint for $\tau_{MM}$
in section \ref{sec:VirMM}. During the proof, we also obtain other non-trivial identities for $A_\la$ in
Lemmas \ref{lem:Phi22M} -- \ref{lem:Phi22M7} and Propositions \ref{prop:Phi21} and \ref{prop:Phi22}. The main result of this subsection
is the following
\begin{thm}\label{thm:MML2}
Let $\Psi(\mu)$ be the function of $\mu$ defined by equation \eqref{eqn:Psi}. Then $\Psi(\mu)=0$
for all strict partitions $\mu$.
\end{thm}
{\bf Proof}:
Note that $\Psi(\mu)$ is skew symmetric with respect to permutations of $\mu$.
So to prove $\Psi(\mu)=0$ for a strict partition $\mu$, it is equivalent to prove
$\Psi(\nu)=0$  where $\nu$ is a permutation of $\mu$. Thus, without loss of generality, we may assume $\mu$ has the standard form  given by equation (\ref{eqn:muStform}).

We first use Corollary \ref{cor:A2} to replace all $A_{2\la}$ in $\Psi(\mu)$ by $A_{\la}$. Then after divided by constant
\[  \frac{(-1)^{r} (1/3)^{(|\mu|+4)/3}}{\prod_{i=1}^p (2k_i-1)!! \prod_{i=1}^q (2m_i-1)!! \prod_{i=1}^r (2n_i+1)!!}, \]
equation $\Psi(\mu)=0$ becomes
\[-\delta_{r+2,q} 3 \Psi_1(\mu) + \delta_{r,q+1}\Psi_2(\mu)=0,\]
where
\begin{align}\label{redu Phi21}
\Psi_1(\mu)=&\frac{9}{2}\sum_{i=1}^q(m_i+1)(6m_i+7)_{[[2]]} A_{\mu+4\epsilon_{i+p}} +\frac{105}{32}A_{(\mu,5,2)}\nonumber\\
&-\frac{1}{16}\sum_{i=1}^q\frac{(6m_i+13)_{[[5]]}}{2m_i+3}A_{\mu+7\epsilon_{i+p}}
\end{align}
and
{\allowdisplaybreaks
\begin{align}\label{redu phi22}
\Psi_2(\mu)
:=&6\sum_{i=1}^p \frac{(6k_i+5)_{[[3]]}}{2k_i+1}
          \bigg( \sum_{j=1 \atop j\neq i}^p  k_j A_{\mu+3\epsilon_{i}+\epsilon_{j}} -\sum_{j=1}^r (n_j+1) A_{\mu+3\epsilon_{i}+\epsilon_{j+p+q}} \bigg)
      \nonumber\\
&
-18\sum_{j=1}^r (n_j+1) \bigg( \sum_{i=1}^q (6m_i+7)_{[[2]]}  A_{\mu+3\epsilon_{i+p}+\epsilon_{j+p+q}}
             +  \sum_{i=1 \atop i\neq j}^r (6n_i+7)_{[[2]]}  A_{\mu+3\epsilon_{i+p+q}+\epsilon_{j+p+q}} \bigg)
\nonumber\\
&+18 \sum_{j=1}^p k_j \bigg( \sum_{i=1}^q(6m_i+7)_{[[2]]}  A_{\mu+3\epsilon_{i+p}+\epsilon_{j}}
               +\sum_{i=1}^r (6n_i+7)_{[[2]]} A_{\mu+3\epsilon_{i+p+q}+\epsilon_{j}} \bigg) \nonumber\\
&+\frac{1}{2}\sum_{i=1}^p\frac{(18k_i^2+45k_i+26)(6k_i+5)_{[[3]]}}{2k_i+1}A_{\mu+4\epsilon_{i}} \nonumber\\
&
-\frac{3}{2}\sum_{i=1}^r(18n_i^2+69n_i+68)(6n_i+7)_{[[2]]}A_{\mu+4\epsilon_{i+p+q}}
 + 6 \sum_{i=1}^p (9k_i+4)A_{(\mu+\epsilon_{i},3)} \nonumber\\
&
-6 \sum_{i=1}^r (9n_i+5)A_{(\mu+\epsilon_{i+p+q},3)}
+3\sum_{i=1}^p \sum_{j=1}^{l(\mu)}(6k_i+1)A_{\mu+\epsilon_{i}+3\epsilon_j}
-3\sum_{i=1}^p(6k_i+1)A_{\mu+\epsilon_{i}} \nonumber\\
&
-3\sum_{i=1}^r\sum_{j=1}^{l(\mu)}(6n_i+5)A_{\mu+\epsilon_{i+p+q}+3\epsilon_j}
+3\sum_{i=1}^r(6n_i+5)A_{\mu+\epsilon_{i+p+q}}
+120A_{(\mu,4)} \nonumber\\
&
-\frac{1}{48}\sum_{i=1}^p\frac{(6k_i+13)_{[[7]]}}{(2k_i+3)(2k_i+1)}A_{(\mu+7\epsilon_{i})}
+\frac{1}{48}\sum_{i=1}^r\frac{(6n_i+17)_{[[7]]}}{(2n_i+5)(2n_i+3)}A_{(\mu+7\epsilon_{i+p+q})} \nonumber\\
&-\frac{15015}{32}A_{(\mu,7)}+\frac{1155}{32}A_{(\mu,6,1)}+\frac{525}{32}A_{(\mu,4,3)}.
\end{align}
}
In calculating $\Psi_1(\mu)$, we have used Lemma $\ref{lem:qi+1}$ to simplify the expression.
In calculating $\Psi_2(\mu)$ we have used equations \eqref{eqn:Ala1} and \eqref{eqn:Ala21} to
remove all terms of the form $(\la,1)$ and $(\la,2,1)$  for some partitions $\la$. Presence of
such terms would make induction process more complicated.

 We will prove $\Psi_1(\mu)=0$  and $\Psi_2(\mu)=0$
 in Propositions \ref{prop:Phi21} and \ref{prop:Phi22} respectively.
 The theorem is thus proved.
$\Box$

\begin{pro}\label{prop:Phi21}
	For $\Psi_1(\mu)$ defined by equation $(\ref{redu Phi21})$, we have
	\[\Psi_1(\mu)=0\]
	for all $\mu$ which is a permutation of a strict partition.
\end{pro}
{\bf Proof}:
Assume $\mu$ is given by equation \eqref{eqn:muStform}.

If $q \neq r+2$, then every term in $\Psi_1(\la)$ is equal to 0 by equation \eqref{eqn:cong}, so we may assume $q = r+2$. We prove this proposition by induction on $p$ and $r$.

If $p=r=0$, then $q=2$. In this case $\Psi_1(\mu)$ is equal to
\begin{eqnarray}
&& \Psi_1((3m_1+1,3m_2+1)) \nonumber \\
&=& \ \frac{9}{2}\Bigg((m_1+1)(6m_1+7)_{[[2]]} A_{(3m_1+5,3m_2+1)}
+(m_2+1)(6m_2+7)_{[[2]]}A_{(3m_1+1,3m_2+5)}\Bigg)   \nonumber \\
&& - \frac{1}{16}\Bigg(\frac{(6m_1+13)_{[[5]]}}{2m_1+3}A_{(3m_1+8,3m_2+1)}
             +\frac{(6m_2+13)_{[[5]]}}{2m_2+3}A_{(3m_1+1,3m_2+8)}\Bigg)  \nonumber \\
&& + \frac{105}{32} \{-A_{(3m_1+1,5)}A_{(3m_2+1,2)}+A_{(3m_1+1,2)}A_{(3m_2+1,5)}\}. \label{eqn:psi1q2}
\end{eqnarray}
In this formula, we have used equation \eqref{eqn:ARec2} to expand $A_{(3m_1+1, 3m_2+1,5,2)}$.
After using equation \eqref{eqn:Al<3} to compute all $A_{\la}$ on the right hand side of this equation and
factoring out $\frac{3}{m_1!m_2!}(\frac{2}{3})^{m_1+m_2+2}$, the first two terms on the right hand side of this equation
together give $-108(m_1-m_2)$ and the last four terms give $108(m_1-m_2)$. Hence we have
\begin{eqnarray} \label{eqn:psi1q2=0}
\Psi_1((3m_1+1,3m_2+1)) &=& 0
\end{eqnarray}
for any non-negative integers $m_1$ and $m_2$.

If $p$ is a positive even integer,  we can apply equation \eqref{eqn:ARec2} to obtain
\begin{align*}
\Psi_1(\mu)
=\sum_{m=2}^p(-1)^mA_{(3k_1,3k_m)}\cdot\Psi_1(\mu^{\{1,m\}}).
\end{align*}

If $p$ is a positive odd integer, we replace each partition $\nu$ appeared in the right hand side of the definition of $\Psi_1(\mu)$ by $(\nu, 0)$ and then use equation \eqref{eqn:ARec2} to obtain
\begin{align*}
\Psi_1(\mu)
=\sum_{m=2}^p(-1)^mA_{(3k_1,3k_m)}\cdot\Psi_1(\mu^{\{1,m\}})+A_{(3k_1,0)}\cdot\Psi_1(\mu^{\{1\}}).
\end{align*}

Since the numbers of parts in $\mu^{\{1\}}$ and $\mu^{\{1, m\}}$ which can be divided by 3 are less than $p$,
the proposition is reduced to the case of $p=0$.

If $r>0$ and $p=0$, we first use recursion formula \eqref{eqn:ARec2} to expand each term in the right hand side of $\Psi_1(\mu)$ in equation \eqref{redu Phi21}. After this expansion, we obtain two terms containing factors $A_{(\mu^{\{1\}},s)}$ with $s=2$ or 5. We then use formula \eqref{eqn:ARec2r} to expand such factors again with respect to the last part. Comparing the result with  equation \eqref{eqn:psi1q2}, we obtain
\begin{align*}
\Psi_1(\mu)
=&\sum_{i=1}^r(-1)^{q+i}A_{(3m_1+1,3n_i+2)}\cdot\Psi_1(\mu^{\{1,q+i\}}) \\
& +\sum_{j=2}^q(-1)^{j}A_{\mu^{\{1,j\}}} \Psi_1((3m_1+1,3m_j+1)).
\end{align*}
Since the second term on the right hand side is $0$ by equation \eqref{eqn:psi1q2=0}, and the numbers of parts in $\mu^{\{1,q+i\}}$ which are equal to 2 mod(3) is $r-1$, the proposition is reduced to the $r=p=0$ case. The proposition is thus  proved. $\Box$

To prove $\Psi_2(\mu)=0$, we will need the following seven lemmas.
\begin{lem}\label{lem:Phi22M}
For $\lambda$ given by equation \eqref{eqn:Stform}, define
\begin{align*}
M_1(\lambda):=\sum_{i=1}^r(-1)^{i}(\frac{2}{3})^{n_i}\frac{1}{n_i!}
\bigg( \sum_{j=1}^p k_j A_{(\lambda+\epsilon_{j})^{\{p+q+i\}}}
-\sum_{j=1 \atop j\neq i}^r (n_j+1) A_{(\lambda+\epsilon_{j+p+q})^{\{p+q+i\}}} \bigg).
\end{align*}
Then $M_1(\lambda)=0$ for all positive partitions $\lambda$.
\end{lem}
{\bf Proof}:
If $r\neq q+2$, every term in $M_1(\lambda)$ is zero by equation \eqref{eqn:cong}. So we may assume $r=q+2$. We will prove this lemma by induction on $p$ and $q$.

The first non-trivial case is $p=q=0$. In this case, $r=2$ and $M_1(\la)$ is equal to
\begin{align*}
M_1((3n_1+2,3n_2+2))
=&\left(\frac{2}{3}\right)^{n_1}\frac{1}{n_1!}\cdot(n_2+1)A_{(3n_2+3)}
-\left(\frac{2}{3}\right)^{n_2}\frac{1}{n_2!}\cdot(n_1+1)A_{(3n_1+3)}\\
=&0,
\end{align*}
where the last equality follows from equation \eqref{eqn:A3k}.

If $p$ is a positive odd integer, then by equation \eqref{eqn:ARec2}, we obtain
\begin{align}\label{eqn:M1_rec}
M_1(\lambda)= & \sum_{i=2}^p(-1)^iA_{(3k_1,3k_i)} M_1(\lambda^{\{1,i\}})
             +\sum_{i=1}^r(-1)^{i} \left( \frac{2}{3} \right)^{n_i}\frac{1}{n_i!} \sum_{j=1 \atop j\neq i}^r (-1)^{\tilde{j}(i)+1+p+q} \nonumber\\
            & \hspace{80pt} \cdot A_{\lambda^{\{1,p+q+i,p+q+j\}}} \left( k_1A_{(3k_1+1,3n_j+2)}-(n_j+1)A_{(3k_1,3n_j+3)} \right),
\end{align}
where $\tilde{j}(i)$ is defined by equation \eqref{eqn:j-tilde}.
By equations \eqref{eqn:Ap2} and \eqref{eqn:Al<3}, the second term on the right hand side of this equation is equal to
\begin{align*}
& (-1)^{p+q+1} \left( \frac{2}{3} \right)^{k_1+1}\frac{1}{k_1!}
\sum_{i,j=1 \atop j\neq i}^r(-1)^{\tilde{j}(i)+i} A_{\lambda^{\{1,p+q+i,p+q+j\}}}\cdot
             \left(\frac{2}{3} \right)^{n_i+n_j}\frac{1}{n_i!n_j!}
=0,
\end{align*}
where the last equality follows from the fact that each summand in the above equation is skew symmetric with respect to $i$ and $j$.
Hence we have
\begin{equation} \label{eqn:M1_rec2}
M_1(\lambda)= \sum_{i=2}^p(-1)^iA_{(3k_1,3k_i)}M_1(\lambda^{\{1,i\}}).
\end{equation}

If $p$ is a positive even integer, we first replace each partition $\nu$ appeared in the definition of $M_1(\lambda)$ by $(\nu,0)$, then apply equation \eqref{eqn:ARec2}. The calculations are similar to the $p$ odd case except that an extra term $A_{(3k_1,0)}M_1(\lambda^{\{1\}})$ should be added to equations \eqref{eqn:M1_rec} and \eqref{eqn:M1_rec2}.

Since the numbers of parts in $\la^{\{1\}}$ and $\la^{\{1,i\}}$ which can be divided by 3 are less than $p$, this lemma can be reduced to the case of $p=0$.

If $q>0$, and $p=0$, we replace each partition $\nu$ appeared in the definition of $M_1(\lambda)$ by $(\nu,0)$ and
expand each term of $M_1(\lambda)$ using equation \eqref{eqn:ARec2}.  Since $A_{(3m_1+1,0)}=0$ , the last part of
$(\nu,0)$ actually does not contribute to the expansion. Hence we have
\begin{align*}
M_1(\lambda)&=\sum_{i=1}^r(-1)^{q+i}A_{(3m_1+1,3n_i+2)}M_1(\lambda^{\{1,i+q\}}).
\end{align*}
Therefore the lemma can be reduced to the $p=q=0$ case and is thus proved.
$\Box$

\begin{lem}\label{lem:Phi22N}
For $\lambda$ given by equation \eqref{eqn:Stform}, define
\begin{align*}
M_2(\lambda):=\sum_{i=1}^r(-1)^{i} \left( \frac{2}{3} \right)^{n_i}\frac{1}{(n_i+1)!}
                \bigg( \sum_{j=1}^p A_{(\lambda+\epsilon_{j})^{\{p+q+i\}}}
                       +\sum_{j=1 \atop j\neq i}^r A_{(\lambda+\epsilon_{j+p+q})^{\{p+q+i\}}} \bigg).
\end{align*}
Then $M_2(\lambda)=0$ for all positive partitions $\lambda$.
\end{lem}
The proof of this lemma is similar to the proof of Lemma \ref{lem:Phi22M}. So we omit it here.

\begin{lem} \label{lem:Phi22M3}
For $\lambda$ given by equation \eqref{eqn:Stform},
define
\[ M_3(\la):=
\sum_{i,j=1 \atop i\neq j}^p (-1)^{j}
(6k_i+5)(6k_i+1)  \frac{(2/3)^{k_j}}{k_j!} A_{(\lambda+3\epsilon_i)^{\{j\}}}.
\]
Then $M_3(\la)=0$ for positive $\la$ with $p>1$ odd.
\end{lem}
{\bf Proof}:
By Corollary \ref{cor:A1&2}, we can factor out $A_{\la_{[2]}}$ from $M_3(\la)$. Hence we
may assume $q=r=0$.

If $p$ is odd, we can use equation \eqref{eqn:ARec2r} to expand $A_{(\lambda+3\epsilon_i)^{\{j\}}}$
 with respect to the part $3k_i +3$ and obtain

\begin{align*}
M_3 (\la)
=&\sum_{i,j=1 \atop i\neq j}^p (-1)^{j}(6k_i+5)(6k_i+1)
\frac{(2/3)^{k_j}}{k_j!} \sum_{m=1 \atop m\neq i,j}^p A_{(3k_i+3,3k_m)}(-1)^{\tilde{i}(j)+\tilde{m}(i,j)}A_{\lambda^{\{i,m,j\}}},
\end{align*}
where $\tilde{i}(j)$ is defined by equation \eqref{eqn:j-tilde}, and
\begin{equation} \label{eqn:m-tilde}
 \tilde{m}(i,j):=m-\delta_{i>m}-\delta_{j<m}
\end{equation}
with $\delta_{i>m}$ being $1$ if $i>m$ and $0$ otherwise.

Computing $A_{(3k_i+3,3k_m)}$ using equation \eqref{eqn:Ap2} and re-grouping terms, we have
\begin{align*}
M_3 (\la)=&\sum_{\{a,b,c\}\subseteq\{1,...,p\} \atop a<b<c}  \frac{(2/3)^{k_a+k_b+k_c+1}}{k_a!k_b!k_c!}A_{\lambda^{\{a,b,c\}}}  \\
& \hspace{80pt} \cdot \sum_{(i,j,m) \in P(a,b,c)}(-1)^{j+\tilde{i}(j)+\tilde{m}(i,j)}
                                    \frac{(6k_i+5)(6k_i+1)(k_i+1-k_m)}{(k_i+1)(k_i+1+k_m)}\\
=&0,
\end{align*}
where the last equality follows from an
elementary identity, whose proof will be given in Lemma \ref{lem:abc2} in the appendix. In the above formula,
$P(a,b,c)$ denotes the set of all permutations of $(a,b,c)$.
$\Box$

\begin{lem} \label{lem:Phi22M4}
For $\lambda$ given by equation \eqref{eqn:Stform},
define
\[ M_4(\la) := \sum_{i,j=1 \atop i\neq j}^p  (-1)^{j} \frac{(2/3)^{k_j}}{k_j!} A_{(\lambda+3\epsilon_i)^{\{j\}}}.
\]
Then $M_4(\la)=0$ for positive $\la$ with $p>1$ odd.
\end{lem}
{\bf Proof}:
By Corollary \ref{cor:A1&2}, we can factor out $A_{\la_{[2]}}$ from $M_4(\la)$. Hence we
may assume $q=r=0$.

If $p$ is odd, we can use equation \eqref{eqn:ARec2r} to expand $A_{(\lambda+3\epsilon_i)^{\{j\}}}$
 with respect to the part $3k_i +3$ and obtain
\begin{align*}
 M_4 (\la)
=&\sum_{i,j=1 \atop i\neq j}^p (-1)^{j}
\frac{(2/3)^{k_j}}{k_j!} \sum_{m=1 \atop m\neq i,j}^pA_{(3k_i+3,3k_m)}(-1)^{\tilde{i}(j)+\tilde{m}(i,j)}A_{\lambda^{\{i,m,j\}}}\\
=&\sum_{\{a,b,c\}\subseteq\{1,...,p\} \atop a<b<c}  \frac{(2/3)^{k_a+k_b+k_c+1}}{k_a!k_b!k_c!}A_{(\lambda^{\{a,b,c\}}} \sum_{(i,j,m)\in P(a,b,c)}(-1)^{j+\tilde{i}(j)+\tilde{m}(i,j)}
\frac{(k_i+1-k_m)}{(k_i+1)(k_i+1+k_m)}\\
=&0,
\end{align*}
where the last equality follows from an
elementary identity, whose proof will be given in Lemma \ref{lem:abc1} in the appendix.
$\Box$

\begin{lem} \label{lem:Phi22M5}
For $\lambda$ given by equation \eqref{eqn:Stform},
define
\[ M_5 (\la):= \sum_{i=1}^p \frac{(-1)^{i} (2/3)^{k_i}}{k_i!} \left(A_{(\lambda,3)^{\{i\}}}-\frac{2}{3}A_{\lambda^{\{i\}}}\right). \]
Then $M_5(\la)=0$ for positive $\la$ with $p$ odd.
\end{lem}
{\bf Proof}:
By Corollary \ref{cor:A1&2}, we can factor out $A_{\la_{[2]}}$ from $M_5(\la)$. Hence we
may assume $q=r=0$.
By equation \eqref{eqn:A3k}, $M_5 (\la)$ can be rewritten as
\[ M_5 (\la)= \sum_{i=1}^p (-1)^{i} A_{(3k_i)} \left(A_{(\lambda,3)^{\{i\}}}-\frac{2}{3}A_{\lambda^{\{i\}}}\right). \]
Since $p$ is odd, we can use Lemma \ref{lem:even-1} for partition $(\la,3)$ to obtain
\begin{align}\label{eqn:lem5_}
\sum_{i=1}^p (-1)^i A_{(3k_i)}A_{(\la,3)^{\{i\}}} = - A_{(3)}A_{\la}.
\end{align}
Since $A_{(3)}= \frac{2}{3}$, we have
\[ M_5 (\la) = \frac{2}{3} \bigg( - A_{\la} + \sum_{i=1}^p (-1)^{i+1} A_{(3k_i)}  A_{\lambda^{\{i\}}} \bigg), \]
which is equal to $0$ by recursion formula \eqref{eqn:ARecOdd0}.
Thus this lemma is proved.
$\Box$

\begin{lem} \label{lem:Phi22M6}
	For $\lambda$ given by equation \eqref{eqn:Stform},
	define
	\[ M_6 (\la):= \sum_{i=1}^p \bigg\{(-1)^{i+1} A_{\la^{\{i\}}} A_{(3k_i+3)}+A_{\la+3\epsilon_i}\bigg\}. \]
	Then $M_6(\la)=0$ for positive partitions $\la$ with $p$ even.
\end{lem}
{\bf Proof}:
By Corollary \ref{cor:A1&2}, we can factor out $A_{\la_{[2]}}$ from $M_6(\la)$. Hence we
may assume $q=r=0$. We prove this lemma by induction on $p$.

The first non-trivial case is $p=2$. In this case $M_6(\la)$ is equal to
\begin{align} \label{eqn:M6ini}
M_6((3k_1,3k_2))=&(A_{(3k_2)}A_{(3k_1+3)}+A_{(3k_1+3,3k_2)})+(-A_{(3k_1)}A_{(3k_2+3)}+A_{(3k_1,3k_2+3)}),
\end{align}
which is zero by straightforward calculations using equations \eqref{eqn:A3k} and \eqref{eqn:Ap2}.

For any even integer $p>2$, we first replace $\la^{\{i\}}$ by $(\la^{\{i\}}, 0)$, then expand
all terms in $M_6(\la)$ with partitions of length bigger than $1$
 using recursion formula \eqref{eqn:ARec2}. We obtain
\begin{align*}
M_6(\la)=&\sum_{j=2}^p (-1)^j A_{(3k_1,3k_j)}M_6(\la^{\{1,j\}})
+\sum_{i=2}^p (-1)^{i+1} A_{(3k_1)}A_{\la^{\{1,i\}}}A_{(3k_i+3)}\\
&+A_{\la^{\{1\}}} A_{(3k_1+3)} +\sum_{i=2}^p (-1)^iA_{\la^{\{1,i\}}}(A_{(3k_1+3,3k_i)}+A_{(3k_1,3k_i+3)}).
\end{align*}
Since $A_{\la^{\{1\}}}=\sum_{i=2}^p(-1)^iA_{(3k_i)}A_{\la^{\{1,i\}}}$ by recursion formula \eqref{eqn:ARecOdd0},
comparing with equation \eqref{eqn:M6ini}, we can see that
\begin{align*}
M_6(\la)=&\sum_{j=2}^p (-1)^j A_{(3k_1,3k_j)}M_6(\la^{\{1,j\}})
+\sum_{i=2}^p (-1)^i A_{\la^{\{1,i\}}}M_6((3k_1,3k_i)).
\end{align*}
This lemma is proved by induction. $\Box$

\begin{lem} \label{lem:Phi22M7}
	For $\lambda$ given by equation \eqref{eqn:Stform},
	define
	\[ M_7 (\la):= \sum_{i=1}^p \bigg\{
	(6k_i^2+6k_i+1)\big( (-1)^{i+1}A_{\la^{\{i\}}} A_{(3k_i+3)}-A_{\la+3\epsilon_i}\big)
	+(-1)^{i+1}A_{(3k_i,3)}A_{\la^{\{i\}}}\bigg\}. \]
	Then $M_7(\la)=0$ for positive $\la$ with $p$ even.
\end{lem}
The proof of this lemma is similar to the proof of Lemma \ref{lem:Phi22M6}. So we omit it here.

Now we are ready to prove the following
\begin{pro}\label{prop:Phi22}
For $\Psi_2(\mu)$ defined by equation $(\ref{redu phi22})$, we have
\[\Psi_2(\mu)=0\]
for all $\mu$ which is a permutation of a strict partition.
\end{pro}
{\bf Proof}:
By equation \eqref{eqn:cong}, the proposition is trivial if $r \neq q+1$.
So we can assume $r=q+1$. We will prove this proposition in four steps.

{\bf Step 1}. We first show that the proposition can be reduced to the case $q=0$
by induction on $q$.

If $q>0$, we apply recursion formula \eqref{eqn:ARec2r} with $j=p+1$ to expand each term $A_\nu$ appeared in the right hand side of equation \eqref{redu phi22}. In case $\nu$  has odd length, we need to replace $\nu$ by $(\nu, 0)$ first. Since in this situation, $\nu_{p+1}$ always has the form $3r+1$ for some non-negative integer $r$ and since $A_{(3r+1,a)} \neq 0$ only if $a=3b+2$ for some integer $b$, many terms in this expansion vanish. In particular, since $A_{(3r+1,0)}=0$, no extra terms would appear in the expansion
when replacing $\nu$ by $(\nu, 0)$. Hence the expansion formulas for $p$ odd case and for $p$ even case are the same.
In case that some $\nu$ has length bigger than $l:=p+q+r$,
$\nu_i$ must be one of the numbers in $\{1, 3, 4, 6, 7\}$ for $i>l$. Hence
$A_{(\nu_{p+1}, \nu_{i})}=0$ for $i>l$. So corresponding terms will not appear in the expansion.
In a summary, after the expansion, we obtain
{\allowdisplaybreaks
\begin{align}
\Psi_2(\mu)=&\sum_{i=1}^r (-1)^{2p+q+i}A_{(3m_1+1,3n_i+i)}\Psi_2(\mu^{\{p+1,p+q+i\}}) \nonumber \\
&+18 \sum_{i=1}^r (-1)^{2p+q+i} \left\{
      (6m_1+7)_{[[2]]} A_{(3m_1+4,3n_i+2)} +(6n_i+7)_{[[2]]} A_{(3m_1+1,3n_i+5)} \right\} \nonumber \\
&   \hspace{30pt} \cdot \bigg( \sum_{j=1}^p k_j A_{(\mu+\epsilon_{j})^{\{p+1,p+q+i\}}}
        - \sum_{j=1 \atop j \neq i}^r(n_j+1) A_{(\mu+\epsilon_{j+p+q})^{\{p+1,p+q+i\}}} \bigg)  \nonumber \\
&+3 \sum_{i=1}^r (-1)^{2p+q+i} \left\{ A_{(3m_1+4,3n_i+2)}+A_{(3m_1+1,3n_i+5)} \right\}  \nonumber \\
& \hspace{30pt}   \cdot \bigg( \sum_{j=1}^p(6k_j+1) A_{(\mu+\epsilon_{j})^{\{p+1,p+q+i\}}}
                               - \sum_{j=1 \atop j \neq i}^r (6n_j+5) A_{(\mu+\epsilon_{p+q+j})^{\{p+1,p+q+i\}}} \bigg).
                               \label{eqn:psi2ind}
\end{align}
}

By equation \eqref{eqn:Al<3}, we have
\begin{align}
&(6m_1+7)_{[[2]]} A_{(3m_1+4,3n_i+2)}
 +(6n_i+7)_{[[2]]}  A_{(3m_1+1,3n_i+5)}  \nonumber \\
& \hspace{120pt} = \frac{72 \cdot (2/3)^{m_1+n_i+2}}{m_1!n_i!}
    - \frac{ 2 \cdot (2/3)^{m_1+n_i+2}}{(m_1+1)!(n_i+1)!}, \label{eqn:Mcoeff}
\end{align}
and
\begin{align}
& A_{(3m_1+4,3n_i+2)}+A_{(3m_1+1,3n_i+5)}
=\frac{ 2 \cdot (2/3)^{m_1+n_i+2}}{(m_1+1)!(n_i+1)!}. \label{eqn:Ncoeff}
\end{align}
Recall that  $M_1(\la)$ and $M_2(\la)$ have been defined in Lemmas \ref{lem:Phi22M} and \ref{lem:Phi22N} respectively.
The first term on the right hand side of equation \eqref{eqn:Mcoeff} combined with corresponding term in
equation \eqref{eqn:psi2ind} produces a factor equal to  $M_1(\mu^{\{p+1\}})$.
The second term on the right hand side of equation \eqref{eqn:Mcoeff} is equal to
the right hand side of equation \eqref{eqn:Ncoeff}. Hence the corresponding terms in equation \eqref{eqn:psi2ind}
can be combined together to produce a factor equal to $M_2(\mu^{\{p+1\}})$.
Putting all these together,  we obtain
\begin{align*}
\Psi_2(\mu)
=&\sum_{i=1}^r (-1)^{2p+q+i}A_{(3m_1+1,3n_i+i)}\Psi_2(\mu^{\{p+1,p+q+i\}})\\
&+\frac{1296(-1)^{2p+q}(2/3)^{m_1+2}}{m_1!} M_1(\mu^{\{p+1\}}) + \frac{6(-1)^{2p+q}(2/3)^{m_1+2}}{(m_1+1)!} M_2(\mu^{\{p+1\}}).
\end{align*}
Since $M_1(\mu^{\{p+1\}})=M_2(\mu^{\{p+1\}})=0$ by Lemmas \ref{lem:Phi22M} and \ref{lem:Phi22N}, this shows that the proposition can be reduced to the $q=0$ case.
Hence for the rest part of the proof for this proposition, we can assume $q=0$, which also implies $r=1$.

{\bf Step 2}. We show $\Psi_2(\mu)=0$ if $p=q=0$ or $p=1$ and $q=0$. For both cases, we have $r=1$.

For the first case, $\mu=(3n+2)$ for some non-negative integer $n$. After using equation~\eqref{eqn:ARecOdd0} to expand $A_{(3n+2,6,1)}$ and $A_{(3n+2,4,3)}$, we have
\begin{align}\label{eqn:psi2_l1}
\Psi_2((3n+2))=&-\left(\frac{3}{2}(18n^2+69n+68)(6n+7)_{[[2]]}+3(6n+5)\right)A_{(3n+6)}
 \nonumber\\
&-6(9n+5)A_{(3n+3,3)} +3(6n+5)A_{(3n+3)}
+\frac{(6n+17)_{[[7]]}}{48(2n+5)(2n+3)}A_{(3n+9)} \nonumber\\
&-\frac{15015}{32}A_{(3n+2,7)}-\frac{1155}{32}A_{(3n+2,1)}A_{(6)}+\frac{525}{32}A_{(3n+2,4)}A_{(3)}
+120A_{(3n+2,4)},
\end{align}
which is zero by straightforward calculations using equations \eqref{eqn:A3k}, \eqref{eqn:Ap2}, and \eqref{eqn:Al<3}.

For the second case, $\mu=(3k,3n+2)$ for some non-negative integers $k$ and $n$. After using equations \eqref{eqn:ARec2} and \eqref{eqn:ARecOdd0} to expand $A_{\nu}$ with $l(\nu)>2$, we have
\begin{align}\label{eqn:psi2_l2}
\Psi_2&((3k,3n+2))=f_1(k,n)+f_2(k,n),
\end{align}
where
\begin{align*}
f_1(k,n):=&A_{(3k+3,3n+3)}\left\{-18(6k+5)(6k+1)(n+1)-3(6n+5)\right\}\\
&+A_{(3k+1,3n+5)}\left\{18(6n+7)(6n+5)k+3(6k+1)\right\}
+6 (9k+4) A_{(3k+1,3n+2)}A_{(3)} \\
&+A_{(3k+4,3n+2)}\left(\frac{3}{2}(6k+5)(6k+1)(18k^2+45k+26)+3(6k+1)\right)\\
&+A_{(3k,3n+6)}\left(-\frac{3}{2}(6n+7)(6n+5)(18n^2+69n+68)-3(6n+5)\right)\\
&- 6(9n+5)\{A_{(3k,3n+3)}A_{(3)}+A_{(3n+3,3)}A_{(3k)}-A_{(3n+3)}A_{(3k,3)}\}\\
&-3(6k+1)A_{(3k+1,3n+2)} +3(6n+5)A_{(3k,3n+3)}
+120A_{(3n+2,4)}A_{(3k)}
\end{align*}
and
\begin{align*}
f_2(k,n):=&-\frac{(6k+13)_{[[7]]}}{48(2k+3)(2k+1)}A_{(3k+7,3n+2)}
+\frac{(6n+17)_{[[7]]}}{48(2n+5)(2n+3)}A_{(3k,3n+9)} \nonumber\\
&-\frac{15015}{32}A_{(3n+2,7)}A_{(3k)}-\frac{1155}{32}A_{(3n+2,1)}A_{(3k,6)}+\frac{525}{32}A_{(3n+2,4)}A_{(3k,3)}.
\end{align*}
A straightforward calculation using  equations \eqref{eqn:A3k}, \eqref{eqn:Ap2},  \eqref{eqn:Al<3} shows that
\[ f_1(k,n)=-f_2(k,n)=\frac{54\cdot(2/3)^{k+n+2}}{k!n!}(36k^2-36kn +6k+18n^2+51n+35).\]
Hence
\begin{align}\label{eqn:ini-Psi2}
\Psi_2&((3k,3n+2))=0
\end{align}
for all non-negative integers $k$ and $n$.

{\bf Step 3}. We prove $\Psi_2 (\mu) =0$ if $q=0$ and $p$ is an odd integer which is bigger than $1$.

Since $r=1$ in this case, we can write $\mu=(\lambda,3n+2)$ where $n$ is a non-negative integer and
 $\la=(3k_1, \cdots, 3k_p)$
for some positive integers $k_1, \cdots, k_p$.

We first use equation $\eqref{eqn:ARec2r}$  to expand each $A_{\nu}$
with respect to $(p+1)$-th part for all partitions $\nu=(\nu_1, \cdots, \nu_{l(\nu)})$
appeared in the right hand side of equation \eqref{redu phi22}.
In case some partitions $\nu$ have odd length, we need replace $\nu$ by $(\nu,0)$ first and then use the equation  $\eqref{eqn:ARec2r}$. For most cases, $\nu_{p+1}$ is equal to $3n+2$ or $3n+5$.
Occasionally, we also have $\nu_{p+1}$ equal to $3n+3$, $3n+6$, or $3n+9$.
If $\nu_{p+1}$ is $3n+2$ or $3n+5$, then $A_{(\nu_{p+1}, \nu_i)}=0$ for all
$\nu_i$ which are not of the form $3r+1$ for some integer $r$. Since
all parts of $\la$ can be divided by $3$,
most terms in the expansion of $A_{\nu}$ vanish.
Moreover, if after the first expansion, we obtain some terms containing
factors like $A_{(\la, s)}$ for some non-negative integer $s$,
we further expand $A_{(\la, s)}$ with respect to the last part $s$
using equation $\eqref{eqn:ARec2r}$.
Actually $s$ can only be $0$, $3$, or $6$ in our case.
In a summary, after such expansions, every term in $\Psi_2 (\mu)$
must contain one of the following three factors:
\[ A_{(\lambda+3\epsilon_i)^{\{j\}}},  \hspace{10pt}
A_{(\lambda,3)^{\{i\}}} , \hspace{10pt}
A_{\lambda^{\{i\}}}
\]
for some $i$ and $j$. Eventually  $\Psi_2(\mu)$ can be written as
\begin{eqnarray}
\Psi_2(\mu)
&=& \sum_{i,j=1,\atop i\neq j}^p (-1)^{p+j} A_{(\lambda+3\epsilon_i)^{\{j\}}} \cdot \omega_{1}(k_i,k_j,n)  +\sum_{i=1}^p (-1)^{p+i} A_{(\lambda,3)^{\{i\}}}\cdot\omega_{2}(k_i,n)
\nonumber \\
&& +\sum_{i=1}^p (-1)^{p+i} A_{\lambda^{\{i\}}} \cdot \{\Psi_2((3k_i,3n+2))- A_{(3)} \omega_{2}(k_i,n)\},
      \label{eqn:Psi_as_psi}
\end{eqnarray}
where
{\allowdisplaybreaks
\begin{align}
\omega_{1}(k_i,k_j,n)
=& \left\{ 18(6k_i+5)(6k_i+1)k_j +3(6k_j+1) \right\} A_{(3k_j+1,3n+2)} \nonumber \\
&-\left\{18(6k_i+5)(6k_i+1)(n+1) + 3(6n+5) \right\}   A_{(3k_j,3n+3)}, \label{eqn:omega1} \\
\omega_{2}(k_i,n)=&6(9k_i+4)A_{(3k_i+1,3n+2)}-6(9n+5)A_{(3k_i,3n+3)}.
\end{align}
}

By straightforward calculations using equations \eqref{eqn:Ap2} and \eqref{eqn:Al<3}, we have
\begin{align}\label{eqn:w1_sim}
\omega_{1}(k_i,k_j,n)=&18(6k_i+5)(6k_i+1)\frac{(2/3)^{k_j+n+1}}{k_j!n!}
+\frac{3 \cdot (2/3)^{k_j+n+1}(6n+7)}{k_j!(n+1)!},
\end{align}
and
\begin{align}\label{eqn:w2_sim}
\omega_{2}(k_i,n)=&\frac{6 \cdot (2/3)^{k_i+n+1}(9n+13)}{k_i!(n+1)!}.
\end{align}

Plugging equations \eqref{eqn:w1_sim}, \eqref{eqn:w2_sim}, \eqref{eqn:ini-Psi2} into equation \eqref{eqn:Psi_as_psi}, we have
\begin{align*}
\Psi_2(\mu)
= \frac{(-1)^{p}}{(n+1)!} \left( \frac{2}{3} \right)^{n+1} \left\{ 18(n+1) M_3(\la)
                                                     +  3(6n+7) M_4(\la)
                                                     +  6(9n+13) M_5(\la) \right\},
\end{align*}
where $M_3(\la), M_4(\la), M_5(\la)$ are defined in
Lemmas  \ref{lem:Phi22M3}, \ref{lem:Phi22M4}, \ref{lem:Phi22M5}, and
they are equal to $0$ by these lemmas.
Hence $\Psi_2(\mu)=0$ if $q=0$ and $p$ is odd.

{\bf Step 4}. We prove $\Psi_2 (\mu) =0$ if $q=0$ and $p$ is a positive even integer.

As in step 3, we expand each term $A_{\nu}$ in the right hand side of equation \eqref{redu phi22}
with respect to $(p+1)$-th part using recursion formula \eqref{eqn:ARec2r}.
In case $\nu$ has odd length, we need replace $\nu$ by $(\nu,0)$ first before expansion.
If after the first expansion, we obtain a term like $A_{(\la, s)}$ with $s=3$ or $6$,
we further expand the term with respect to the $(p+1)$-th part after replacing
$(\la, s)$ by $(\la, s,0)$. The calculation is similar to step 3 except that in case
we obtain terms like $A_{\la}$ and $A_{\la + 3 \epsilon_i}$ after first expansion, we do not expand these terms further.
Eventually, after the expansion, every term contains one of the following five factors:
\[A_{(\lambda+3\epsilon_i)^{\{j\}}},\ \ A_{\lambda^{\{i\}}},\ \ A_{(\la,3)^{\{i\}}}, \ \ A_{\lambda+3\epsilon_i},\ \ A_{\lambda}\]
for some $i$ and $j$, and $\Psi_2(\mu)$ can be written as
\begin{eqnarray}
\Psi_2(\mu)
&=& \sum_{i,j=1,\atop i\neq j}^p (-1)^{p+j} A_{(\lambda+3\epsilon_i)^{\{j\}}} \cdot \omega_{1}(k_i,k_j,n)
+ A_{\lambda} \cdot \Psi_2((3n+2)) \nonumber \\
&& +\sum_{i=1}^p (-1)^{p+i} A_{\lambda^{\{i\}}} \cdot \xi_{1}(k_i,n)
+\sum_{i=1}^p A_{\lambda+3\epsilon_i} \cdot \xi_{2}(k_i,n) \nonumber \\
&&+\sum_{i=1}^p (-1)^{p+i} A_{(\lambda,3)^{\{i\}}} \cdot \xi_{3}(k_i,n),
\label{eqn:psitow1234}
\end{eqnarray}
where $\omega_1$ is defined by equation \eqref{eqn:omega1}, and
{\allowdisplaybreaks
\begin{align*}
\xi_{1}(k_i,n)=&\Psi_2((3k_i,3n+2))- A_{(3)} \omega_{2}(k_i,n)
-6(9n+5)A_{(3n+3)}A_{(3k_i,3)}\\
&+ \left( 6(9n+5)A_{(3n+3,3)} - 120A_{(3n+2,4)} + \frac{15015}{32} A_{(3n+2,7)} \right) \cdot A_{(3k_i)},\\
\xi_{2}(k_i,n)
=&- \left\{ 18(6k_i+5)(6k_i+1)(n+1) + 3(6n+5) \right\} A_{(3n+3)},\\
\xi_{3}(k_i,n)
=&6(9k_i+4)A_{(3k_i+1,3n+2)}
-6(9n+5)A_{(3k_i,3n+3)}.
\end{align*}
}

By equation \eqref{eqn:w1_sim}, we have
\[\omega_{1}(k_i,k_j,n)=\left( 18(6k_i+5)(6k_i+1)\frac{(2/3)^{n+1}}{n!}
+\frac{3 \cdot (2/3)^{n+1}(6n+7)}{(n+1)!}\right)\cdot A_{(3k_j)}.\]
When considering the contribution of $\omega_{1}(k_i,k_j,n)$ to $\Psi_2(\mu)$ in equation \eqref{eqn:psitow1234}, for each $i$, we apply Lemma \ref{lem:even-1} for $(\la+3\epsilon_i)$ and obtain
\[\sum_{j=1  \atop j\neq i}^p (-1)^{j+1} A_{(3k_j)} A_{(\la+3\epsilon_i)^{\{j\}}}
= (-1)^{i} A_{(3k_i+3)} A_{\la^{\{i\}}}. \]
This formula will be used to simplify the first term on the right hand side of equation \eqref{eqn:psitow1234}.

To compute the contribution of $\xi_{1}(k_i,n)$ in equation \eqref{eqn:psitow1234}, we use equation \eqref{eqn:w2_sim} to rewrite $\omega_2(k_i,n)$ as
\begin{align*}
\omega_{2}(k_i,n)
=\frac{6 \cdot (2/3)^{n+1}(9n+13)}{(n+1)!}\cdot A_{(3k_i)}.
\end{align*}
Hence, by equation \eqref{eqn:ini-Psi2},  we have
\begin{equation} \label{eqn:tw3}
\xi_{1}(k_i,n) = -6(9n+5)A_{(3n+3)}A_{(3k_i,3)} + c(n) A_{(3k_i)},
\end{equation}
where $c(n)$ is a function only depending on $n$.
By Lemma \ref{lem:even-1},
\[\sum_{i=1}^p (-1)^{p+i} A_{\la^{\{i\}}} \cdot c(n) A_{(3k_i)}=0.\]
Hence the second term on the right hand side of equation \eqref{eqn:tw3} does not contribute to the calculation of
$\Psi_2(\mu)$ in equation \eqref{eqn:psitow1234}.

By straightforward calculations,
\begin{align*}
\xi_{3}(k_i,n)=\frac{(2/3)^{k_i+n+1}}{k_i!(n+1)!}\cdot6(9n+13)
=6(9n+13)\frac{(2/3)^{n+1}}{(n+1)!}\cdot A_{(3k_i)}.
\end{align*}
When calculating its contribution to $\Psi_2(\mu)$ in equation \eqref{eqn:psitow1234}, we use the following
formula
\[\sum_{i=1}^p (-1)^{i} A_{(3k_i)}A_{(\lambda,3)^{\{i\}}}= - \sum_{i=1}^p(-1)^{i}A_{(3k_i,3)}A_{\la^{\{i\}}}
,\]
which follows from the fact that both sides of this equation are equal to $A_{(3)}A_{\la}-A_{(\la,3)}$
by equation \eqref{eqn:ARecOdd0} applied to $A_{(\la, 3)}$ and equation\eqref{eqn:ARec2r} applied to $A_{(\la,3,0)}$ with
$j=p+1$.

Putting above considerations together, we obtain the following formula from equation \eqref{eqn:psitow1234}:
\begin{align}
\Psi_2(\mu)\label{eqn:psi2_peven}
=&\frac{(2/3)^{n+1}}{(n+1)!}\cdot\bigg\{\sum_{i=1}^p (-1)^{i+1} A_{\la^{\{i\}}} A_{(3k_i+3)}
\cdot \left\{18(6k_i+5)(6k_i+1)(n+1)+3(6n+7) \right\} \nonumber \\
& \hspace{60pt}
-\sum_{i=1}^p A_{\lambda+3\epsilon_i} \cdot \left\{ 18(6k_i+5)(6k_i+1)(n+1)+3(6n+5) \right\} \nonumber \\
& \hspace{60pt} -108(n+1)\sum_{i=1}^p (-1)^{i} A_{\lambda^{\{i\}}}A_{(3k_i,3)} \bigg\}.
\end{align}
Recall that $M_6$ and $M_7$ have been defined in Lemmas \ref{lem:Phi22M6} and \ref{lem:Phi22M7} respectively. Writing $(6n+7)=6(n+1)+1$ and $(6n+5)=6(n+1)-1$ in equation \eqref{eqn:psi2_peven}, and grouping all  terms
containing factor $(n+1)$ together, we obtain
\begin{align*}
\Psi_2(\mu)=&3\frac{(2/3)^{n+1}}{(n+1)!}\cdot M_6 (\la)+108\frac{(2/3)^{n+1}}{n!}\cdot M_7(\la).
\end{align*}
Since $ M_6 (\la)= M_7 (\la)=0$ by Lemmas \ref{lem:Phi22M6} and \ref{lem:Phi22M7}, this proposition is proved.
 $\Box$

This completes the proof of Theorem \ref{thm:MML2}.

\appendix
\vspace{30pt}
\hspace{160pt} {\bf \Large Appendix}

\section{Proof of Corollary \ref{cor:pr-nota}}
\label{sec:prdr}

In this section we give a proof of Corollary \ref{cor:pr-nota}.

We first show equation \eqref{eqn:multi pr} holds for all partitions. By skew symmetry of
Q-polynomials, both sides of equation \eqref{eqn:multi pr}
are skew symmetric with respect to permutations of $\la$ if $\la$ is weakly positive.
If $\la$ has at least two positive equal parts, then equation \eqref{eqn:multi pr} holds since both sides are zero.
By Lemma \ref{lem:multi pr}, equation \eqref{eqn:multi pr} holds for all positive partitions.
We then prove equation \eqref{eqn:multi pr} by induction on the number of parts of $\la$ which are zero.
Through permutation of a zero part and positive parts, we can assume
 $\la=(\mu, 0)$ where the number of zero parts in $\mu$ is strictly less than the number of zero parts
 in $\la$. In this case,  equation \eqref{eqn:multi pr} has
the form
\[ rt_rQ_{(\mu, 0)}=\sum_{i=1}^{l(\mu)}Q_{(\mu+r\epsilon_i, 0)} + Q_{(\mu, r)} + \frac{1}{2}Q_{(\mu, 0, r)}
+\sum_{k=1}^{r-1}\frac{(-1)^{r-k}}{4}Q_{(\mu,0,k,r-k)}. \]
The two middle terms on the right hand side of this equation adds up to $\frac{1}{2} Q_{(\mu, r)}$.
The $0$ parts appeared in partitions in all other terms in the above equation can be simply removed.
The result is exactly the same as equation \eqref{eqn:multi pr} applied to
$\mu$ which holds by induction hypothesis.


From the proof of equation \eqref{eqn:dQ} outlined in \cite{Mac}  page 265-266, we can see that
this equation holds for a larger system of functions $Q_{\la}$
which are defined for all $\la \in \mathbb{Z}^l$.
In particular, it must hold when $\la$ is a partition.
 This completes the proof of Corollary \ref{cor:pr-nota}.
$\Box$

\section{Some elementary identities}

In this section we prove some elementary identities which are needed in the proof of combinatorial identities
in sections \ref{sec:Phi=0} and \ref{sec:Psi=0}.

\begin{lem}\label{lem:abc1}
	For any positive integers $a < b <c$  and non-negative integers $k_a, k_b, k_c$ indexed by $a,b,c$, we have
	\begin{align*}
	\sum_{(i,j,m) \in P(a,b,c)}(-1)^{j+\tilde{i}(j)+\tilde{m}(i,j)}
	\frac{(k_i+1-k_m)}{(k_i+1)(k_i+1+k_m)}=0,
	\end{align*}
where $\tilde{j}(i)$ has been defined by equation \eqref{eqn:j-tilde}, $\tilde{m}(i,j)$ has been defined
by equation \eqref{eqn:m-tilde}, and $P(a,b,c)$ denotes the set of all permutations of $(a,b,c)$.
\end{lem}
{\bf Proof}:
Note that
 \[j+\tilde{i}(j)+\tilde{m}(i,j)=i+j+m-\delta_{i<j}-\delta_{i>m}-\delta_{j<m}.\]
Using the identity
\[\frac{(k_i+1-k_m)}{(k_i+1)(k_i+1+k_m)}
-\frac{(k_m+1-k_i)}{(k_m+1)(k_m+1+k_i)}
=\frac{1}{k_m+1}-\frac{1}{k_i+1},\]
we have
\begin{align*}
& \sum_{(i,j,m) \in P(a,b,c)} (-1)^{j+\tilde{i}(j)+\tilde{m}(i,j)}
\frac{(k_i+1-k_m)}{(k_i+1)(k_i+1+k_m)}\\
=&(-1)^{a+b+c}\left\{\left(\frac{1}{k_c+1}-\frac{1}{k_a+1}\right)
      +\left(\frac{1}{k_b+1}-\frac{1}{k_c+1}\right)
      +\left(\frac{1}{k_a+1}-\frac{1}{k_b+1}\right)\right\}\\
=&0.
\end{align*}
$\Box$

\begin{lem}\label{lem:abc2}
	For any positive integers $a<b<c$ and non-negative integers $k_a, k_b, k_c$ indexed by $a,b,c$, we have
	\begin{align*}
	\sum_{(i,j,m) \in P(a,b,c)}(-1)^{j+\tilde{i}(j)+\tilde{m}(i,j)}D_{i,m}=0,
	\end{align*}
	where  $D_{i,m}:=\frac{(6k_i+5)(6k_i+1)(k_i+1-k_m)}{(k_i+1)(k_i+1+k_m)}$.
\end{lem}
{\bf Proof}:
Let
$E_{i,m}:=\frac{36k_ik_m+36k_i+36k_m+41}{k_i+1}$. Then
\begin{align*}
D_{i,m}-D_{m,i}=E_{m,i}-E_{i,m}.
\end{align*}
Since $E_{i,j}-E_{i,m}=36(k_j-k_m)$, we have
\begin{align*}
 & \sum_{(i,j,m) \in P(a,b,c)} (-1)^{j+\tilde{i}(j)+\tilde{m}(i,j)}D_{i,m}  \\
= & (-1)^{a+b+c+1}\left\{(E_{a,c}-E_{a,b})+(E_{b,a}-E_{b,c})+(E_{c,b}-E_{c,a})\right\} = 0. 
\end{align*}
 $\Box$


\vspace{30pt} \noindent
Xiaobo Liu \\
School of Mathematical Sciences \& \\
Beijing International Center for Mathematical Research, \\
Peking University, Beijing, China. \\
Email: {\it xbliu@math.pku.edu.cn}
\ \\ \ \\
Chenglang Yang \\
Beijing International Center for Mathematical Research, \\
Peking University, Beijing, China. \\
Email: {\it yangcl@pku.edu.cn}
\end{document}